\documentclass[11pt]{article}
\usepackage{amsfonts,amsmath}
\usepackage{latexsym}
\usepackage{amsthm}
\usepackage{amscd}
\usepackage{commath}
\usepackage{epsfig}
\usepackage{amssymb}
\usepackage{tikz}
\usepackage{graphicx}
\usepackage{caption}
\usepackage{enumitem}
\usepackage{mathtools}
\usepackage[toc,page]{appendix}
\newcommand{\tree}{\begin{tikzpicture}
		\draw[line width=0.06cm, black] (0,0) -- (0,1)  {};
		\draw[line width=0.06cm, black] (0,1) -- (1,2)  {};
		\draw[line width=0.06cm, black] (0,1) -- (0,2)  {};
		\draw[line width=0.06cm, black] (0,1) -- (-1,2)  {};
		\draw[line width=0.06cm, black] (-1,2) -- (-1.8,3)  {};
		\draw[line width=0.06cm, black] (-1,2) -- (-0.6,3)  {};
		\draw[line width=0.06cm, black] (1,2) -- (1,3)  {};
		\draw[line width=0.06cm, black] (1,2) -- (1.7,3)  {};
		\draw[line width=0.05cm, purple] (-1,2) -- (-2.2,3)  {};
		\draw[line width=0.05cm, purple] (-2.2,3) -- (-3,4)  {};
		\draw[line width=0.05cm, purple] (-2.2,3) -- (-2.4,4)  {};
		\draw[line width=0.05cm, blue] (-1,2) -- (-1.4,3)  {};
		\draw[line width=0.05cm, blue] (-1,2) -- (-0.95,3)  {};
		\draw[line width=0.05cm, blue] (-0.95,3) -- (-1.3,4)  {};
		\draw[line width=0.05cm, blue] (-0.95,3) -- (-0.8,4)  {};
		\draw[line width=0.05cm, blue] (-1.3,4) -- (-1.9,5)  {};
		\draw[line width=0.05cm, blue] (-1.3,4) -- (-1.3,5)  {};
		\draw[line width=0.05cm, blue] (-1.3,4) -- (-0.8,5)  {};
		\draw[line width=0.05cm, brown] (-1,2) -- (-0.1,3)  {};
		\draw[line width=0.05cm, brown] (-1,2) -- (0.4,3)  {};
		\draw[line width=0.05cm, brown] (0.4,3) -- (0,4)  {};
		\draw[line width=0.05cm, brown] (0.4,3) -- (0.6,4)  {};
		\node [fill=black, inner sep=2.5 pt,label=182: $i$] at (-1,2) {} ;
		\node [fill=black, inner sep=2.5 pt,label=left: $\phi_2(i)$] at (0,1) {} ;
		\node [fill=black, inner sep=2.5 pt,label=182: $i_0$] at (0,0) {} ;
\end{tikzpicture}}

\newcommand{\treeleaf}{\begin{tikzpicture}
		
		\draw[line width=0.05cm, black] (0,1) -- (1,2)  {};
		\draw[line width=0.05cm, black] (0,1) -- (0,2)  {};
		\draw[line width=0.05cm, black] (0,1) -- (-1,2)  {};
		\draw[line width=0.05cm, black] (-1,2) -- (-1.8,3)  {};
		\draw[line width=0.05cm, black] (-1,2) -- (-0.6,3)  {};
		\draw[line width=0.05cm, black] (1,2) -- (1,3)  {};
		\draw[line width=0.05cm, black] (1,2) -- (1.7,3)  {};
		\draw [line width=0.1 cm, olive] plot [smooth cycle] coordinates {(-1.8,3) (-2.4,6) (-1,6)};
		\draw [line width=0.04 cm, red] plot [smooth cycle] coordinates {(-0.6,3) (-0.9,4.5)  (-0.3,4.5)};
		\draw [line width=0.04 cm, red] plot [smooth cycle] coordinates {(1,3) (0.4,5.5) (1.4,5.5)};
		\draw [line width=0.04 cm, red] plot [smooth cycle] coordinates {(1.7,3) (1.7,5) (2.4,5)};
		\node[label=$T_1$] at (-1.83,4) {};
		\node[label=$T_2$] at (-0.6,3.65) {};
		\node[label=$T_3$] at (0.95,4.2) {};
		\node[label=$T_4$] at (1.95,4) {};
		\node [fill=black, inner sep=2.5 pt,label=182: $i_0$] at (0,0) {} ;
		\draw[line width=0.05cm, black] (0,0) -- (0,1)  {};
		\node [fill=black, inner sep=2.5 pt,label=182: $i_1$] at (-1.8,3) {} ;
		\node [fill=black, inner sep=2.5 pt,label=182: $i_2$] at (-0.6,3) {} ;
		\node [fill=black, inner sep=2.5 pt,label=182: $i_3$] at (1,3) {} ;
		\node [fill=black, inner sep=2.5 pt,label=0: $i_4$] at (1.7,3) {} ;
		\node [fill=black, inner sep=1.5 pt] at (0,2) {} ;
		\node [fill=black, inner sep=1.5 pt] at (-1,2) {} ;
		\node [fill=black, inner sep=1.5 pt] at (1,2) {} ;
\end{tikzpicture}}

\newcommand{\treespine}{\begin{tikzpicture}
		\draw[line width=0.07cm, black] (0,0) -- (0,1)  {};
		\draw[line width=0.06cm, black] (0,1) -- (1,2)  {};
		\draw[line width=0.06cm, black] (0,1) -- (-1,2)  {};
		\draw[line width=0.06cm, black] (-1,2) -- (-1.8,3)  {};
		\draw[line width=0.06cm, black] (-1,2) -- (-0.6,3)  {};
		\draw[line width=0.06cm, black] (1,2) -- (0.7,3)  {};
		\draw[line width=0.05cm, black] (1,2) -- (1.7,3)  {};
		\draw [line width=0.08 cm, olive] plot [smooth cycle] coordinates {(-1.8,3) (-2,7) (-1.6,7)};
		\draw [line width=0.04 cm, red] plot [smooth cycle] coordinates {(-0.6,3) (-0.8,7)  (-0.4,7)};
		\draw [line width=0.04 cm, red] plot [smooth cycle] coordinates {(0.7,3) (0.5,7) (0.9,7)};
		\draw [line width=0.04 cm, red] plot [smooth cycle] coordinates {(1.7,3) (1.5,7) (1.9,7)};
		\draw [line width=0.04 cm, blue] plot [smooth cycle] coordinates {(-1,2) (-1.4,4.5) (-0.9,4.5) (-1,3)};
		\draw [line width=0.04 cm, blue] plot [smooth cycle] coordinates {(-1,2) (-0.5,2.6) (-0.25,3.3) (-0.15,3.3) (-0.2,2.5)};
		\draw [line width=0.04 cm, blue] plot [smooth cycle] coordinates {(1,2) (0.5,2.6) (0.3,3.3) (0.1,3.3) (0.2,2.5)};
		\draw [line width=0.04 cm, blue] plot [smooth cycle] coordinates {(1,2) (1,2.5) (1.1,4) (1.4,4)};
		\draw [line width=0.04 cm, blue] plot [smooth cycle] coordinates {(1,2) (1.7,2.7) (3,5) (2.5,3)};
		\draw [line width=0.04 cm, blue] plot [smooth cycle] coordinates {(0,1) (-0.2,2.2) (0.3,2.1) (0.05,1.2) };
		\draw [line width=0.04 cm, blue] plot [smooth cycle] coordinates {(0,1) (1.7,1.9) (2.4,2.2) (1.6,1.5)};
		\node [fill=black, inner sep=2 pt] at (-1.8,3) {} ;
		\node [fill=black, inner sep=2 pt] at (-0.6,3) {} ;
		\node [fill=black, inner sep=2 pt] at (0.7,3) {} ;
		\node [fill=black, inner sep=2 pt] at (1.7,3) {} ;
		\node [fill=black, inner sep=2 pt] at (-1,2) {} ;
		\node [fill=black, inner sep=2 pt] at (1,2) {} ;
		\node [fill=black, inner sep=2 pt] at (0,1) {} ;
		\node [fill=black, inner sep=2.5 pt,label=182: $i_0$] at (0,0) {} ;
\end{tikzpicture}}

\usepackage{float}

\usepackage{color,soul}

\newcommand{\bb}{\begin{equation}}
	\newcommand{\ee}{\end{equation}}

\newcommand{\cT}{\mathcal T}
\newcommand{\cB}{\mathcal B}
\newcommand{\cF}{\mathcal F}

\newcommand{\W}{\mathbb{W}}
\newcommand{\D}{\mathbb{D}}

\definecolor{dblue}{rgb}{.61,.61,1}
\definecolor{lightblue}{rgb}{.61,.61,1}
\definecolor{altblue}{rgb}{.61,.61,1}

\tikzset{cross/.style={cross out, draw=black, minimum size=2*(#1-\pgflinewidth), inner sep=0pt, outer sep=0pt},
	cross/.default={1pt}}

\newtheorem{thm}{Theorem}[section]
\newtheorem{prop}[thm]{Proposition}
\newtheorem{lem}[thm]{Lemma}
\newtheorem{cor}[thm]{Corollary}
\newtheorem*{thm*}{Theorem}
\theoremstyle{remark}
\newtheorem{remark}{Remark}[section]

\theoremstyle{definition}

\def\IB{\relax\hbox{$\inbar\kern-.3em{\rm B}$}}
\def\IC{\relax\hbox{$\inbar\kern-.3em{\rm C}$}}
\def\ID{\relax\hbox{$\inbar\kern-.3em{\rm D}$}}
\def\IE{\relax\hbox{$\inbar\kern-.3em{\rm E}$}}
\def\IF{\relax\hbox{$\inbar\kern-.3em{\rm F}$}}
\def\IG{\relax\hbox{$\inbar\kern-.3em{\rm G}$}}
\def\IGa{\relax\hbox{${\rm I}\kern-.18em\mathcal{T}$}}
\def\IH{\relax{\rm I\kern-.18em H}}
\def\IK{\relax{\rm I\kern-.18em K}}
\def\IL{\relax{\rm I\kern-.18em L}}
\def\IP{\relax{\rm I\kern-.18em P}}
\def\IR{\relax{\rm I\kern-.18em R}}
\def\IZ{\relax{\rm Z\kern-.5em Z}}

\begin{document}
	\tikzset{
		position label/.style={
			below = 3pt,
			text height = 1.5ex,
			text depth = 1ex
		},
		brace/.style={
			decoration={brace, mirror},
			decorate
		}
	}

	\thispagestyle{empty}
	\quad
	
	\vspace{2cm}
	\begin{center}
		
		\textbf{\huge Local limits of one-sided trees}
		
		\vspace{1.5cm}
		
		{\Large Bergfinnur Durhuus ~~~~~~~~ Meltem Ünel}
		\vspace{0.5cm}
		
		{\it Department of Mathematical Sciences, Copenhagen University\\
			Universitetsparken 5, DK-2100 Copenhagen {\O}, Denmark}

		\vspace{0.5cm} {\sf durhuus@math.ku.dk, meltem@math.ku.dk}

		\vspace{1.5cm}
		
		{\large\textbf{Abstract}}\end{center} A finite \emph{one-sided tree} of height $h$ is defined as a rooted planar tree obtained by grafting branches on one side, say the right, of a spine, i.e. a linear path of length $h$ starting at the root, such that the resulting tree has no simple path starting at the root of length greater than $h$.  We consider the distribution $\tau_N$ on the set of one-sided trees $T$ of fixed size $N$, such that the weight of $T$  is proportional to $e^{-\mu h(T)}$, where $\mu$ is a real constant and $h(T)$ denotes the height of $T$. We show that, for $N$ large, $\tau_N$ has a weak limit as a probability measure  supported on infinite one-sided trees. The dependence of the limit  measure $\tau$ on $\mu$ shows a transition at $\mu_0=-\ln 2$ from a single spine phase for $\mu\leq \mu_0$ to a multi-spine phase for $\mu> \mu_0$. Correspondingly, there is a transition in the volume growth rate of balls around the root as a function of radius from linear growth for $\mu<\mu_0$, to quadratic growth at $\mu=\mu_0$, and to qubic growth for $\mu> \mu_0$.

	\vspace{3cm}
	
	\newpage

	\section{Introduction}
	
	Aspects of various families of random trees, such as the classical Bienaymé-Galton-Watson (BGW) trees, regular trees, and planar labelled trees among many others, have been intensively studied in the literature \cite{drmota2009random} for a variety of reasons, including their close connection to branching processes \cite{athreyaney1972branching} and in some cases because of their usefulness in studying statistical behaviour of other graph ensembles via various bijective correspondences, see e.g. \cite{schaeffer1998conj}. Thus, enumeration of such classes of trees and the critical behaviour of the corresponding generating functions is a classical subject \cite{flajolet2009analytic} providing in many cases a basis for the construction of local and scaling limits, see \cite{abraham2015introduction,le2012scaling} for non-exhaustive reviews. Most work in this area concerns so-called simply generated trees, where the weight of individual (finite) trees is a local function of the vertex degrees or, more precisely, can be written as a product of weight functions depending on the individual vertex degrees. In previous work \cite{meltem2021alpha, meltem2021height}, we considered some particular cases of rooted planar trees where the weight function is non-local due to a dependence on the height of individual trees. In the case of an exponential dependence on height, the weight of a tree $T$ is given by
	\bb \label{weight}
	w(T)=	e^{-\mu h(T)} g^{-|T|}\,,
	\ee
	where $|T|$ denotes the size i.e. the number of edges and $h(T)$ the height of $T$, while $g$ and $\mu$ are fixed real parameters. Since $2|T|$ equals the sum of all vertex degrees in $T$, the second factor on the right-hand side is local while the height-dependent factor is not, except for $\mu=0$ which corresponds to the uniform distribution of trees of fixed size. For this family of weights, a nontrivial dependence on the strength $\mu$ of the height-coupling was found \cite{ meltem2021height}, exhibiting a transition at $\mu=0$ between different asymptotic volume behaviours of balls around the root at large radius. Specifically, this was accomplished through constructing local limits at large size of the corresponding ensembles. In the present paper we consider in a similar vein a restricted class of rooted planar trees, so-called \textit{one-sided trees} admitting a spine of maximal height and branches grafted on its right side, and show that a similar behaviour occurs but with the important difference that the critical coupling is shifted from $\mu=0$ to $\mu =-\ln 2$. In particular, uniformly distributed trees exhibit quite different behaviours in the two instances, as will be seen. 
	
	Apart from providing interesting examples of random trees with non-local weight functions whose local limits can be constructed and analysed in some detail, a partial motivation for investigating trees with weights of the form \eqref{weight} also arises in the study of certain loop models on two-dimensional \textit{causal triangulations} \cite{durhuus2021critical} through a  bijective correspondence between such causal triangulations of the disk or the cylinder and rooted planar trees \cite{malyshev2001two,durhuus2010spectral}. Indeed, it turns out that the so-called dilute loop model, in the limit of small coupling, via this correspondence is equivalent to the unrestricted random planar tree model with weight function \eqref{weight} corresponding to $\mu= -\ln 2$. From quite a different perspective, the two-dimensional causal triangulation model with so-called curvature dependent coupling has been studied by other techniques  in \cite{di2000integrable}, defined in terms of a subclass of cylindrical causal triangulations that are in bijection with a set of triangulations of the rectangle subject to so-called \textit{staircase boundary conditions} \cite{di2000integrable}. Without entering into a more detailed explanation of this concept, it suffices in the present context to note that the tree correspondence in this case yields a bijection between the latter triangulations of fixed size and one-sided trees of fixed size, and hence providing additional motivation for initiating a study of their properties. 
	
	
	The paper is organized as follows. In subsection 2.1, we review some notions associated with rooted planar trees and define the set $\Omega$ of one-sided trees as a subset of the set of rooted planar trees $\cT$. A convenient metric $\text{dist}(\cdot,\cdot)$ defining a topology and an associated Borel $\sigma$-algebra on $\cT$ is introduced, yielding an appropriate setting for discussing local (or weak) limits. Moreover, the notion of \textit{grafting} and the \textit{spine} map are briefly discussed. In subsection 2.2, the analytic structure of the generating functions $Y_m$ for the numbers $B_{m,N}$ of one-sided trees of fixed height $m$ and of size $N$ is determined and used to derive closed expressions for those numbers, of importance for the analysis in section \ref{sec:5}. The finite size probability measures $\tau_N ^{(\mu)}, N \in \mathbb{N}$, are then defined by restricting \eqref{weight} to one-sided trees of fixed size $N$ and normalizing. The relevant normalisation factors are denoted by $W_N^{(\mu)}$ and their asymptotic behaviour for large $N$ is of crucial importance for the determination of the local limit.

	Section \ref{sec:3} covers the case $\mu<\mu_0$, using here and in the following the shorthand 
	$$
	\mu_0= -\ln 2\,.
	$$
	\noindent A rather simple analysis of the analytic properties of the generating function $W^{(\mu)}(g)$ for the partition functions $W^{(\mu)}_N$ allows a  determination of their asymptotic behaviour  in subsection \ref{sec:3.1}, while in section \ref{sec:3.2} lower bounds on ball volumes in $\cT$ are established as a prerequisite for determining the local limit $\tau^{(\mu)}$ of $(\tau^{(\mu)}_N)_{N \in \mathbb{N}}$. Finally, a linear increase of ball volumes around the root as a function of radius is shown.
	
	In section \ref{sec:4}, the case $\mu=\mu_0$ is considered. The proof of  Proposition \ref{prop:ZNln2} on the logarithmic singularity of the partition function $\W^{(\mu_0)} (g)$ is deferred to a appendix. Moreover, the identification of the local limit with a BGW measure as well as the quadratic volume growth of balls around the root is shown.
	
	Section \ref{sec:5} is devoted to the case $\mu> \mu_0$. Results from \cite{meltem2021height} are used to determine the asymptotic behaviour of $W^{(\mu)} _N$ in subsection \ref{sec:5.1}, while existence of the local limit $\nu^{(\mu)}$ is established  in section \ref{sec:5.2} on the basis of ball volume extimates in $\cT$. Basic properties of $\nu^{(\mu)}$ are investigated in subsection 5.3; in particular, the spine is in this case shown to be multi-ended and its governing measure is shown to be equal to the one found in \cite{meltem2021height} up to a shift of the coupling $\mu$ by $\ln 2$, see Proposition \ref{pushforward}. Finally, the volume growth of the balls around the root is shown to be cubic in Corollary \ref{expected}.

	\section{Preliminaries}\label{sec:2}
	
	\subsection{Spaces of trees}\label{sec:2.1}
	
	In the following we consider various spaces of \textit{planar rooted trees}, to which we simply refer as \textit{trees}. We assume that the reader is familiar with the basic notions associated with trees and refer to section 2 of \cite{meltem2021height} for more details, where similar notation is used. In particular, we denote the vertex set of a tree $T$ by $V(T)$ or just $V$ and similarly the edge set is called $E(T)$ or $E$. The \emph{root vertex} $i_0$ of the trees considered is always assumed to have degree $1$ and the degree $\sigma(v)$ of any vertex $v$ is assumed to be finite. The space of all such planar trees will be denoted by $\cT$, with the subsets $\cT_{\rm fin}$ and $\cT_\infty$ consisting of finite and infinite trees, respectively. For $N\in\mathbb N$ we denote by ${\mathcal T}_N$ the set of planar trees $T$ of size $N$, i.e. $|T|:=|E(T)|=N$, where $|A|$ is used to designate the cardinality of any set $A$. Using the graph distance, the set of vertices in $T$ at distance $r$ from the root will be called $D_r(T)$, and we note that because of planarity, each $D_r$ is naturally ordered from left to right, when using the standard orientation of the plane. In particular, $D_0(T)=\{i_0\}$ and $D_1(T)=\{i_1\}$, where $i_1$ is the unique neighbour of the root, and the edge $\{i_0, i_1\}$ will be called the \emph{root edge}. Given a vertex $v\in D_r(T)$, where $r\geq 1$, it has a unique \emph{ancestor} (or \emph{parent}) in $D_{r-1}(T)$, that we denote by $\phi_r(v)$, thus defining an order preserving surjective map $\phi_r: D_r(T)\to D_{r-1}(T)$. The elements of $\phi_r^{-1}(u)$ are called the \emph{offspring} of $u\in D_{r-1}(T)$. 
	
	The height $h(T)$ of a tree $T\in\cT_{\rm fin}$ is defined as  the maximal distance of any vertex in $T$ from the root, or
	$$
	h(T) = {\rm max}\{r\mid D_r(T)\neq \emptyset\}\,.
	$$
	If $T\in \cT_\infty$ we set $h(T)=\infty$.
	
	Considering simple paths in a tree $T\in\cT$ starting at the root, there is obviously a unique left-most path $\omega_T$ of maximal length, in the sense that the $k$th vertex in $\omega_T$ is the left-most vertex in $D_k(T)$ for each finite $k=0,1,...,|\omega_T|$, where $|\omega_T|$ denotes the length of $\omega_T$, i.e. the number of edges in $\omega_T$ (defined to be $\infty$ if $\omega_T$ is infinite). We say that $T$ is a \emph{one-sided tree}  if $|\omega_T|=h(T)$. The subsets of $\cT_N$, $\cT_{\infty}$, $\cT_{\rm fin}$ and $\cT$ consisting of one-sided trees will be denoted by $\Omega_N$, $\Omega_{\infty}$, $\Omega_{\rm fin}$ and $\Omega$ respectively. Moreover, we use an upper index in parenthesis on a set of trees to indicate the height, e.g. $\Omega^{(m)}$ stands for the subset of $\Omega$ consisting of one-sided trees of height $m$ and $\Omega_N^{(m)}$ is the set of one-sided trees of size $N$ and height $m$.
	
	In order to discuss local limits of sequences of measures on $\cT$ we equip it with a natural metric as follows.   
	Given $r\in\mathbb N$, let $B_r(T)$ denote the ball of radius $r$ around the root $i_0$ defined as the subtree of $T$ spanned by the vertices at distance at most $r$ from $i_0$.
	For $T,T'\in\cT$ we then set 
	$$
	\text{dist}(T,T') = \inf\{\frac{1}{r} \mid r\in\mathbb N,\, B_r(T)=B_r(T')\} \,. 
	$$
	
	It is then easy to see that dist is a metric on $\cT$, in fact an ultrametric. We shall denote by $\mathcal{B}_{a}(T_0)$  the ball of radius $a>0$ around $T_0\in\cT$, 
	$$
	\cB_a(T_0) := \{ T \in \cT \mid \text{dist}(T,T_0) \leq a \} \, .
	$$
	
	The measures on $\cT$ discussed in the following are all Borel measures, i.e. they are defined on the Borel $\sigma$-algebra $\cF$ generated by the open sets. By definition, a sequence $(\tau_N)_{N\in\mathbb N}$ of probability measures converges weakly to a probability measure $\tau$ on $\cT$, if 
	$$
	\int_\cT F\,d\tau_n\;\to\; \int_\cT F\,d\tau \quad \mbox{as $N\to\infty$}
	$$
	for all real valued bounded continuous functions $F$ on $\cT$. This requirement is equivalent to the statement that
	\bb\label{convcond}
	\tau_N(\cB)\;\to\; \tau(\cB)\quad\mbox{as $N\to\infty$}
	\ee
	for any ball $\cB$ in $\cT$, as further detailed in the following remark listing some basic properties of the metric space $\cT$. They are easily verifiable and will be used repeatedly in the subsequent discussion, see e.g.  \cite{chassaing2006local,durhuus2003probabilistic} for more details. 
	\begin{remark}\label{remcT}

		i)\; Any ball in $\cT$ is both open and closed and any two balls are either disjoint or one is contained in the other. Since $\cT_{\rm fin}$ is a countable dense subset of $\cT$, it follows immediately by use of Theorem 2.2 in \cite{billingsley2013convergence} that a sequence of probability measures $(\tau_N)_{N\in\mathbb N}$ converges weakly to a probability measure $\tau$ on $\cT$, if \eqref{convcond} holds for all balls $ \cB$. 
		
		ii)\;  If $T\neq T'$, then $\text{dist}(T,T') = \frac {1}{r}$, where $r\geq 1$ is the radius of the largest ball around their roots shared by $T$ and $T'$, and in this case we have 
		$$
		\cB_{\frac 1r}(T) = \cB_{\frac 1r}(T') = \cB_{\frac 1r}(T_0)\,,\quad\mbox{where $T_0=B_r(T)$ and hence $h(T_0)=r$}.
		$$ 
		\noindent Note that if $r> h(T)$, then $\cB_{\frac 1r}(T)= \{T \}$. It follows that all other balls can be written in the form $\cB_{\frac 1r} (T)$ where $r=h(T)$ and they generate $\mathcal{F}$.	
		
		iii)\; $\cT$ is a complete metric space. It is not compact, but subsets of the form 
		\bb \nonumber  \label{Compact}
		C = \bigcap _{r=1} ^\infty \{ T \in \mathcal T \big| ~ |D_r (T)|  \leq K_r \}\,,
		\ee
		where $(K_r)_{r\in\mathbb N}$ is any sequence of positive numbers, are compact. 
		
		iv)\; $\Omega_\infty$ is a closed subspace of $\cT$. 
	\end{remark}
	Given $T\in\cT_\infty$, we say that a vertex $i$ of $T$ is of \emph{infinite type} if it has infinitely many descendants, and otherwise it is of \emph{finite type}. Clearly, the vertices of infinite type span a subtree of $T$ with the same root and root edge and with no leaves, i.e. no vertices of degree 1. We call this subtree the  \emph{spine} or \emph{backbone} of $T$ and denote it by $\chi(T)$.  The mapping $\chi:\cT_\infty\to\cT_\infty$ will be called the \emph{spine map}. As shown in \cite{meltem2021height}, $\chi$ is a Borel measurable map and its image $\cT_s := \chi(\cT_\infty)$ is a closed subset of $\cT$.
	
	We will make use of the notion of \emph{grafting} a tree $T_1\in \cT$ onto another tree $T_0\in\cT$ in the following. To be specific, let $i$ be some vertex in $T_0$ (different from the root) of degree $\sigma(i)$, and note that the edges of $T_0$ incident on $i$ divide a sufficiently small disc around $i$ into $\sigma(i)$ sectors that can be ordered clockwise around $i$ and which we shall denote by $S_{(i,n)},\,1\leq n\leq\sigma_i$.  The grafted tree $T := gr(T_0;(i,n);T_1)$ is then defined 
	by identifying the oriented root edge $(i_0,i_1)$ of $T_1$ with the edge $(\phi_{r_1}(i),i)$ of $T_0$ and drawing the remaining part of $T_1$ in the $n$'th sector $S_{i,n}$ of the plane around $i$, see Fig.1. Clearly, $T_0$ and $T_1$ can be considered as subtrees of $T$, although $T_1$ and $T$ have different roots unless $i=i_1$. We say that $T$ is obtained by \emph{grafting $T_1$ onto $T_0$ at $(i,n)$}. In case $i$ is a leaf, there is only one sector $S_{i,1}$ and we say that $T_1$ is grafted onto $T_0$ at $i$.
	
	It is easily seen that, for fixed $T_0$ and pairs $(i_1,n_1),\dots,(i_K,n_K)$ labelling different vertex sectors,  successive grafting of trees $T_1,\dots, T_K$ at $(i_1,n_1),\dots, (i_K,n_K)$, respectively, is well defined and independent of the order of grafting. We denote the so obtained tree by $gr(T_0;(i_1,n_1),\dots,(i_K,n_K);T_1,\dots,T_K)$. 
	
	\begin{remark} \label{remarkgraft}
		It is easy to verify that the mapping
		$$
		G: (T_1,\dots,T_K)\to gr(T_0;(i_1,n_1),\dots,(i_K,n_K);T_1,\dots,T_K)
		$$
		maps $\cT^K$ homeomorphically onto a closed subset of $\cT$, and if $T_0\in \cT_{\rm fin}$, the image is also open. Moreover, if $T_0$ is finite and $i_1,\dots, i_K$ denote the vertices (leaves) at maximal height $r:=h(T_0)$, then $G$ maps $\cT^K$ homeomorphically onto $\cB_{\frac 1r}(T_0)$. We refer to section 2 of \cite{meltem2021height} for more detailed arguments.
	\end{remark}
	
	\begin{figure}[ht]
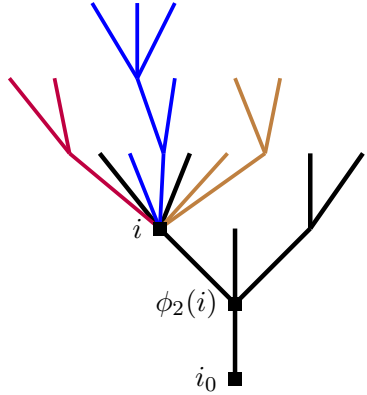

		\centering
		\begin{align*}
			\raisebox{0cm}{\tree} 
		\end{align*}
		\captionof{figure}{A finite tree $T_0$ (edges in black) with root $i_0$ and subtrees  \textcolor{purple}{$T_1$}, \textcolor{blue}{$T_2$}, \textcolor{brown}{$T_3$} grafted in the three sectors at $i$. }
		\label{figure:tree}
	\end{figure}
	
	\subsection{Generating functions}\label{sec:2.2}

	We define the generating function for the numbers $A_{m,N}$ of trees of height at most $m$ and size $N$ by 
	\bb\label{partfuncfinite}
	X_m(g) = \sum_{h(T)\leq m} g^{|T|} = \sum_{N=1}^\infty A_{m,N} g^N
	\ee
	and, with  $A_N$ denoting the total number of trees of size $N$, we set 
	\bb\label{defX}
	X(g) =  \sum_{T\in \cT_{\rm fin}} g^{|T|} = \sum_{N=1}^\infty A_{N} g^N\,.
	\ee
	It is well known that $X_m$ and $X$ satisfy 
	\bb\label{eq:XXm}
	X_m(g) = \frac{g}{1-X_{m-1}(g)}\,,\quad m\geq 2\,,\qquad  X(g) = \frac{g}{1-X(g)}\,.
	\ee
	From the latter equation follows that the series \eqref{defX} converges for $|g|<\frac 14$ and the sum is given by
	\bb\label{solX}
	X(g) = \frac{1-\sqrt{1-4g}}{2}\,.
	\ee
	In particular, this implies that $A_N$ equals the $(N-1)$'th Catalan number:
	\bb\label{asympCatalan} \nonumber
	A_N = C_{N-1} := \frac{(2N-2)!}{N!(N-1)!}\,. 
	\ee
	Using that $X_1(g)=g$, the former equation in \eqref{eq:XXm} can be solved recursively, yielding (see e.g. \cite{drmota2009random} for details)
	\bb\label{eq:sol1}
	X_m(g) = 2g\frac{(1+\sqrt{1-4g})^m -(1-\sqrt{1-4g})^m}{(1+\sqrt{1-4g})^{m+1} - (1-\sqrt{1-4g})^{m+1}} \,,\quad m\geq 1\,,
	\ee
	which can also be written as
	\bb\label{eq:sol2}
	X_m(g) = \sqrt g\frac{U_{m-1}(\frac{1}{2\sqrt g})}{U_m(\frac{1}{2\sqrt g})}\,,
	\ee
	where $U_m$ denotes  the $m$'th Chebychev polynomial of the second kind. 
	
	Given $T\in\Omega^{(m)}\,, m\geq 2$, let $i_2$ denote the \underline{leftmost} descendant of $i_1$ and let $T'$ denote the subtree of $T$ spanned by $(i_1, i_2)$ and all the descendants of $i_2$. Then $T'\in\Omega^{(m-1)}$ with root equal to $i_1$, while the tree obtained by removing $T'$ (except vertex $i_1$) from $T$ belongs to $\cup_{k=1}^m\cT^{(k)}$. Clearly, this decomposition gives a bijective correspondence between trees in $\Omega^{(m)}$ and $(\cup_{k=1}^m \cT^{(k)})\times \Omega^{(m-1)}$. Setting 
	$$
	Y_m(g) = \sum_{T\in\Omega^{(m)}} g^{|T|} := \sum_{N=1}^\infty B_{m,N} g^N\,,
	$$
	this implies
	\bb\label{eq:recursionY}
	Y_m(g) =  X_m(g)Y_{m-1}(g)\,,\quad m\geq 2\,.  
	\ee
	Using  $Y_1(g)=X_1(g)=g$ and \eqref{eq:sol2} we obtain
	\bb\label{Y1}
	Y_m(g) = \prod_{k=1}^m X_k(g) = \frac{g^{\frac m2}}{ U_m(\frac{1}{2\sqrt g})}\,.
	\ee
	Recalling that $U_m$ is a polynomial of degree m with the same parity as $m$ it follows that $Y_m$ is a rational function of $g$ of the form
	\bb\label{Y2} \nonumber
	Y_m(g) = \frac{g^m}{P_m(g)}\,,
	\ee
	where $P_m$ is a polynomial of degree $\lfloor\frac m2\rfloor$ with non-vanishing constant term. From \eqref{Y1} it is easy to determine the poles and residues of $\frac{1}{P_m}$. Indeed, for $m\geq 2$, the poles are simple and located at the nonvanishing roots of $U_m$ that are known to be given by
	$$ 
	g_{m,k}= \frac 14\Big(1+\tan^2\frac{\pi k}{m+1}\Big)\,,\quad k=1,\dots,\big\lfloor\frac m2\big\rfloor\,,
	$$   
	and the corresponding residues are found to be (see e.g. \cite{meltem2021height} for a similar calculation) 
	\bb\label{eq:res2} \nonumber
	r_{m,k} =  \frac{(-1)^k}{m+1}\,g_{m,k}^{\frac{1-m}{2}}\,\tan^2\frac{\pi k}{m+1}\,,\quad k=1,\dots,\big\lfloor\frac m2\big\rfloor\,.
	\ee
	Hence, from
	\bb\label{Y3} \nonumber
	Y_m(g) = g^m\sum_{k=1}^{\lfloor\frac m2\rfloor} \frac{r_{m,k}}{g-g_{m,k}}\,,
	\ee
	we obtain by expanding the pole terms on the right-hand side that, for $N\geq m\geq 2$, 
	\bb\label{BmN}
	B_{m,N} = \frac{2^{-(m-1)}}{m+1} \sum_{k=1}^{\lfloor\frac m2\rfloor} (-1)^{k+1}\tan^2\frac{\pi k}{m+1}\Big(1+\tan^2\frac{\pi k}{m+1}\Big)^{-N+\frac{m-1}{2}}\cdot 4^N\,.
	\ee
	For later use, this formula should be compared to the expression found in \cite{meltem2021height} for the Taylor coefficients $A_{m,N}$ of $X_m$,
	\bb\label{AmN}
	A_{m,N} = \frac{1}{m+1} \sum_{k=1}^{\lfloor\frac m2\rfloor} \tan^2\frac{\pi k}{m+1}\Big(1+\tan^2\frac{\pi k}{m+1}\Big)^{-N}\cdot 4^N\,.
	\ee   
	
	The goal of the subsequent discussion is to study the one-parameter family of probability measures $\tau^{(\mu)}$ obtained as local limits of the finite size measures $\tau^{(\mu)} _N$ concentrated on $\Omega_N$ and defined by
	\bb\label{eq:tauN} \nonumber
	\tau_N^{(\mu)}(T) = \frac{e^{-\mu h(T)}}{W^{(\mu)}_N}\,,\quad\mbox{for $T\in\Omega_N$}\,,
	\ee 
	where the normalisation factor (partition function) $W^{(\mu)}_N$ is given by 
	\bb\label{eq:part}
	{W^{(\mu)}_N} = \sum_{m=1}^\infty e^{-\mu m}B_{m,N}\,.
	\ee 
	We shall consider these as measures on the space $\cT$ and obtain the local limits as measures supported on $\Omega_\infty$. 
	
	For later reference we note that the corresponding construction for the case of unrestricted planar trees was carried out in \cite{meltem2021height}. More precisely, defining the probability measure $\nu^{(\mu)}_N$ on $\cT_N$ by 
	\bb\label{eq:vuN}
	\nu_N^{(\mu)}(T) = \frac{e^{-\mu h(T)}}{Z^{(\mu)}_N}\,,\quad\mbox{for $T\in\cT_N$}\,,
	\ee 
	where 
	\bb\label{eq:partZ}
	{Z^{(\mu)}_N} = \sum_{m=1}^\infty e^{-\mu m}A_{m,N}\,,
	\ee 
	the existence of the limit $\lim_{N\to\infty}\nu^{(\mu)}_N :=\nu^{(\mu)}$ was established and some of its properties uncovered. The particular case $\mu=0$, concerning uniformly distributed trees, has been thoroughly investigated in the literature (see e.g. \cite{abraham2015introduction} and references therein). In particular, it is known that the limit $\tau^{(0)}$, called the Uniform Infinite Planar Tree (UIPT), is a single-spine tree whose branches are critical and i.i.d. as explained in more detail in section \ref{sec:3.3}.  For $\mu>0$, on the other hand, the measure $\nu^{(\mu)}$ is supported on multi-spine trees and induces a Poisson type measure $\tilde\nu^{(\mu)}$ on spine trees, i.e. on $\cT_s$, whose explicit characterisation is given in \eqref{nutilde} below. As shown in section \ref{sec:4}, the measure $\tau^{(\mu)}$ is supported on multi-spine trees for $\mu>\mu_0$ and induces the measure $\tilde\nu^{(\mu+\ln 2)}$ on $\cT_s$. In particular, the uniform infinite planar one-sided tree $\tau^{(0)}$ is a multi-spine tree, contrary to $\nu^{(0)}$.
	
	\section{The case $\mu< \mu_0$}\label{sec:3}
	
	\subsection{Partition Function}\label{sec:3.1}
	
	In order to investigate the asymptotic behaviour of $W_N^{(\mu)}$ for $N$ large we consider, for fixed $\mu$, the generating function 
	\bb\label{def:W}
	W^{(\mu)}(g) = \sum_{m=1}^\infty e^{-\mu m} Y_m(g) = \sum_{m=1}^\infty \sum_{N=1}^\infty B_{m,N}e^{-\mu m}g^N\,,
	\ee
	for which the following holds.
	\begin{prop}\label{prop:gc}
		The radius of convergence $g_c(\mu)$ for the power series defining $W^{(\mu)}$ is given by
		\bb\label{gc} \nonumber
		g_c(\mu) = \begin{cases} \frac 14 \qquad~~~~~~~~ \mbox{if $\mu>\mu_0$} \\ e^{\mu}(1- e^{\mu}) \quad\mbox{if $\mu \leq\mu_0$}\,.\end{cases}
		\ee 
	\end{prop}
	\begin{proof}
		Since the power series \eqref{partfuncfinite} is divergent for $g>g_{m,1} := g_m$, the same holds for $Y_m$ by \eqref{Y1}, and since $g_m\to \frac 14$ as $m\to\infty$, it follows that the series defining $W^{(\mu)}$ diverges for $|g|>\frac 14$. On the other hand,  using \eqref{eq:XXm} we have by Lemma 3.1 in \cite{meltem2021height} that 
		\bb\label{Yconverge}
		Y_m(g) X(g)^{-m} = \prod_{l=1}^m \frac {X_l(g)}{X(g)} = (1-X(g)) \prod_{l=1}^{m-1}\frac{1-X(g)}{1-X_l(g)}\to f(g)\quad \mbox{as $m\to\infty$}
		\ee
		for all $g\in\mathbb D$, where 
		$$\mathbb D := \{g \in \mathbb C \mid |g| < \frac 14 \}$$
		and where $f$ is an analytic function on $\mathbb D$ with no zeroes. Since $|X(g)|\leq X(|g|)$, this implies that $W^{(\mu)}(g)$ is finite provided $e^{-\mu} X(|g|) < 1$, which by \eqref{solX} is seen to be fulfilled for $|g|<\frac 14$ if $\mu>\mu_0$, while for $\mu\leq \mu_0$ it is equivalent to 
		$$
		|g| < e^{\mu}(1- e^{\mu})\,.
		$$
		In the latter case, it also follows from \eqref{Yconverge} that the series \eqref{def:W} diverges for $|g|>e^\mu(1-e^\mu)$. This completes the proof.
	\end{proof}	
	
	\begin{thm}
		For fixed $\mu >\mu_0$, there exists  $b >g_c(\mu)$ such that $W^{(\mu)}(g)$ is analytic in
		\bb \label{puncdisk} \nonumber
		\{g\in\mathbb C\mid \abs{g}< b , g \neq g_c(\mu) \}  \,,
		\ee
		and has a simple pole at $g_c(\mu)$.
		\label{thmZ}
	\end{thm}
	
	\begin{proof}
		Using \eqref{Yconverge}, we have
		\bb\label{Ymterm}
		e^{-\mu m} Y_{m}(g) =  \big(e^{-\mu}X(g)\big)^{m} (1-X(g)) \prod_{l=1} ^{m-1} \frac{1-X(g)}{1-X_l(g)}\,.
		\ee
		Defining 
		\bb\label{defc} \nonumber
		c(g) = \frac{g}{(1-X(g))^2}
		\ee
		we have by Lemma 3.4 in \cite{meltem2021height} that \eqref{Ymterm} can be rewritten as 
		\bb\label{Ymterm2}
		e^{-\mu m} Y_{m}(g) = f(g) \big(e^{-\mu}X(g)\big)^{m}  + h_m(g)\,,
		\ee
		where $h_m$ is analytic in $\mathbb D$ and fulfills 
		\begin{align*}
			|h_m(g)| &= \big| e^{-\mu}X(g)\big| ^{m} \big| (1-X(g))\prod_{l=1}^{m-1} \frac{1-X(g)}{1-X_l(g)}- f(g) \big| \\
			&\leq \mbox{\rm cst}\cdot \big(e^{-\mu}X(|g|) c(|g|)\big)^m\quad\mbox{for $|g|\leq b$ and $m\in\mathbb N$}\,, 
		\end{align*}
		for any fixed  $b<\frac 14$. It is easily seen, that $c(|g|))<1$ for $|g|<\frac 14$ and in particular $c(g_c(\mu)) <1$. Hence, $b>g_c(\mu)$ can be chosen such that $e^{-\mu}X(|g|)c(|g|)<1$ for $|g|\leq b$.
		It then follows that  $\sum_{m=1}^\infty h_m(g)$ converges to an analytic function $h(g)$  for $|g|<b$ and, by summing over $m$ in \eqref{Ymterm2}, we conclude that
		$$
		W^{(\mu)}(g) =  \frac{f(g)}{1-e^{-\mu}X(g)} + h(g) 
		$$
		is analytic for $|g|<b$ except at $g=g_c(\mu)$, which is a simple zero of $1-e^{-\mu}X(g)$.  This completes the proof.  
	\end{proof}
	
	\begin{cor}\label{cor:Zasymp1}
		There exists $d>0$ such that
		\bb
		W^{(\mu)}_N =  r\cdot g_c(\mu)^{-(N+1)} \big(1+ O(e^{-d N})\big)
		\label{Zasymp1}
		\ee
		for $N$ large, where $r$ is the residue of $-W^{(\mu)}(g)$ at $g=g_c(\mu)$.
	\end{cor}
	
	\begin{proof}
		By Theorem \ref{thmZ} we may write 
		$$
		W^{(\mu)}(g) = \frac{r}{g_c(\mu)-g} + \tilde h(g)\,,
		$$
		where $\tilde h$ is analytic in a disc centred at $0$ of radius $b>g_c(k)$. Expanding the pole term as a geometric series in $\frac{g}{g_c(\mu)}$ then yields the dominant term in \eqref{Zasymp1}, while the subdominant part arises from the Taylor coefficients of $\tilde h$.
	\end{proof}

	\subsection{Lower bounds on ball volumes and the local limit}\label{sec:3.2}

	To prove that the sequence $(\tau^{(\mu)} _N)_{N\in\mathbb N}$ given by \eqref{eq:vuN} has a weak limit on $\cT$, we proceed by first establishing lower bounds on ball volumes that will allow us to prove tightness of the measures $\tau^{(\mu)}_N$ and subsequently to identify the limit.

	\begin{lem} \label{lem:ballbound1}
		Let $\mu<\mu_0$ and let $T_0\in\Omega_{\rm fin}$ have height $r$ with $K$ vertices in $D_r(T_0)$. For each $M\in\mathbb N$ there exists $d>0$ such that
		\bb\label{ballbound1}
		\tau^{(\mu)} _N (\cB_{\frac 1r}(T_0)) \geq    e^{-\mu(r-1)} g_c(\mu)^{|T_0|-K}   \Big(\sum_{S=1}^M C_{S-1}g_c(\mu)^{S}\Big)^{K-1} \big(1+O(e^{-d N})\big)\,.
		\ee
	\end{lem}		
	\begin{proof}
		Given $T_0 $ as stated, let $i_1, \dots , i_K$ denote the vertices at maximal height $r=h(T_0)$, ordered from left to right, such that $i_1$ is the lefttmost vertex in $D_r(T_0)$. Recall that any tree $T$ in $\mathcal{B} _{\frac{1}{r}}(T_0)$ can be obtained by grafting $K$ trees $T_1, \dots, T_K \in \cT$ onto $T_0$ at $i_1,\dots, i_K$, respectively, which we shall refer to as the branches of $T$, see Fig.2. If $|T|=N$, we then have   
		\bb\label{eq:sumN}
		N=  |T_0|+|T_1|+\dots+|T_K| -K\,,
		\ee
		and since the measure $\tau_N ^{(\mu)}$ is supported on one-sided trees, we can assume $T\in\Omega_N$ and hence that $T_1$ belongs to $\Omega$ and has maximal height among $T_1,\dots, T_K$ satisfying 
		$$
		h(T) = h(T_1) + r-1\,.
		$$
		Given $M\in\mathbb N$ and imposing the constraints $|T_i|\leq M$ for $i\neq 1$ and $h(T_1)> M$, it follows that $T_1$ is the unique branch of maximal size and, setting $|T_i|=N_i$, it holds that 
		\bb\label{ballboundprel}
		\tau^{(\mu)} _N (\cB_{\frac 1r}(T_0)) \geq  e^{-\mu(r-1)} \sum_{\substack{ N_i \leq M , i>1\\ N_1+...+N_K= N-|T_0|+K }} (W^{(\mu)}_N) ^{-1} W_{M,N_{1}} \prod_{i =2}^{K} C_{N_i-1}  \,,
		\ee
		for $N$ large enough, where 
		$$
		W_{M,N} := \sum\limits_{m=M+1}^\infty e^{-m\mu} B_{m,N}\,.
		$$
		Clearly, the generating function for the $W_{M,N}$ for fixed $M$ equals $W^{(\mu)}(g)$ minus the function
		
		\noindent $\sum\limits_{m=1}^{M-1} e^{-m \mu } Y_m(g)$, which is analytic in $\D$. As a consequence, $W_{N,M}$ has the same asymptotic form \eqref{Zasymp1} for $N$ large as $W^{(\mu)}_N$. For fixed values of $N_i,\,i > 1$,  it follows that the corresponding term in \eqref{ballboundprel} fulfills 
		$$
		(W^{(\mu)}_N) ^{-1} W_{M,N_{1}} = g_c(\mu)^{N-N_1}\big(1+O(e^{-d N})\big)\,.
		$$
		Inserting this into \eqref{ballboundprel} and using \eqref{eq:sumN}, we obtain \eqref{ballbound1}.
	\end{proof}
	
	\begin{figure}[ht]
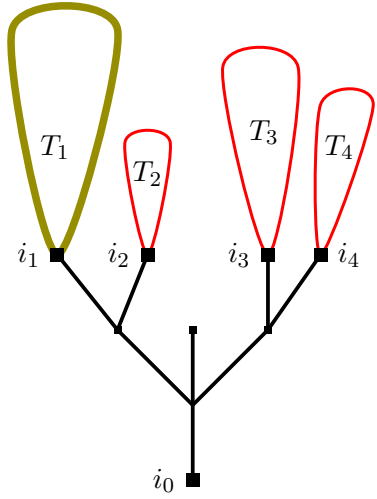

		\centering
		\begin{align*}
			\raisebox{0cm}{\treeleaf} 
		\end{align*}
		\captionof{figure}{Structure of a tree $T$ in $\cB _{\frac{1}{3}}(T_0) \cap \Omega$ around $T_0 \in \Omega^{(3)}$ with root $i_0$ and leaves $(i_1,i_2,i_3,i_4)$ at height $h(T_0)=3$, obtained by grafting a one-sided tree $T_1 \in \Omega$ (in green) at the leftmost vertex $i_1$ and $T_2,T_3,T_4 \in \cT$ (in red) at $i_2,i_3,i_4$, respectively. Note that since $T \in \Omega$, $T_1$ has the maximal height among $T_1,T_2,T_3,T_4$.}
		\label{figure:treeleaf}
	\end{figure}
	
	Let us denote the large-$N$ limit of the right-hand side of \eqref{ballbound1} by $\Lambda(T_0,M)$, i.e.
	\begin{align} \label{defLambdaM} \nonumber
		\Lambda(T_0,M)  =   e^{-\mu(r-1)}  g_c(\mu)^{|T_0| -K} \big(\sum_{S=1}^M C_{S-1} g_c(\mu)^S\Big) ^{K-1}\,, 
	\end{align}
	and define 
	\bb\label{defLambda}
	\Lambda(T_0) := \lim_{M\to\infty} \Lambda(T_0,M) =   e^{-\mu(r-1)}  g_c(\mu)^{|T_0| -K} X(g_c(\mu))^{K-1}\,. 
	\ee
	
	\begin{lem}\label{lem:addupto} For any $r\in\mathbb N$, it holds that 
		\bb\label{addupto}
		\sum_{T_0 \in \Omega^{(r)}} \Lambda(T_0)  \; =\; 1\,.
		\ee
	\end{lem}
	
	\begin{proof}
		We use an inductive argument. For $r=1$ the statement trivially holds, 
		so let $r \geq 2$ be arbitrary and assume \eqref{addupto} holds for $r-1$. For given $T_0\in \Omega^r$, let $T_0^\prime= B_{r-1}(T)$ and denote by $K'$ the number of vertices in $T_0^\prime$ at height $r-1$. Obviously, $T' _0 \in \Omega$ and for fixed such $T^\prime _0$ the number of trees $T_0\in\Omega^r$ with  $K$ vertices at height $K$ equals $\binom{K+K^\prime-2}{K^\prime-1}$. Hence, we have 
		\begin{align} \nonumber
			\sum_{\substack{T_0 : B_{r-1}(T_0) = T_0 ^\prime\\ h(T_0)=r}} \Lambda(T_0) 
			&=  e^{-\mu(r-1)} \cdot g_c(\mu) ^{|T_0^\prime|} \sum_{K \geq 1} \binom{K+K^\prime-2}{K^\prime-1} X(g_c(\mu))^{K-1} \\ \nonumber
			&= e^{-\mu(r-1)} \cdot g_c(\mu)^{|T_0 ^\prime|} \big(1-X(g_c(\mu)\big) ^{-K^\prime} \\
			&= e^{-\mu(r-2)} \cdot g_c(\mu)^{|T_0 ^\prime|-K^\prime} X(g_c(\mu))^{K^\prime -1}\,,  \label{identical}
		\end{align}
		where the second equality follows by using the identity 
		\bb\label{combid}
		\sum_{K=R}^\infty \binom{K}{R}x^K = \frac{x^R}{(1-x)^{R+1}}\,,
		\ee
		and in the final step \eqref{eq:XXm} and  the fact that
		\bb\label{eq:gc}
		e^{-\mu}X(g_c(\mu))=1
		\ee 
		were used. Since the last expression in \eqref{identical} equals $\Lambda(T_0^\prime)$, this completes the proof.
	\end{proof}	
	
	We can now establish the existence of the limit measure and provide a characterisation of it as follows.

	\begin{thm}\label{thm:limit1}
		For each $\mu<\mu_0$ the sequence  $(\tau^{(\mu)}_N)$ is weakly convergent to a Borel probability measure $\tau^{(\mu)}$ on $\cT$ supported on $\Omega_\infty$, that is characterized by 
		\bb\label{limit1}
		\tau ^{(\mu)}(\cB_{\frac 1r}(T_0)) = \Lambda(T_0) = e^{-\mu(r-1)}\cdot g_c(\mu)^{|T_0| -K} X(g_c(\mu))^{K-1}\,,
		\ee
		for any  tree $T_0\in\Omega_{\rm fin}$ of height $r$, where $K=|D_r(T_0)|$.
	\end{thm}
	
	\begin{proof}
		First we note that the sequence  $(\tau_N^{(\mu)})$ is tight as a consequence of Lemmas \ref{lem:ballbound1} and \ref{lem:addupto}. The proof of this fact is essentially identical to that of Corollary 3.6 in \cite{meltem2021height} to which we refer for further details.
		It follows that $(\tau^{(\mu)} _N)$ has a weakly convergent subsequence $(\tau^{(\mu)} _{N_i})$ converging to a probability measure $\tau^{(\mu)} $ on $\cT$. We shall show that the limit $\tau^{(\mu)}$ is independent of the subsequence and hence that $(\tau^{(\mu)} _N)$ is convergent. Since the balls in $\cT$ have empty boundary, we  have by Theorem 2.1  in \cite{billingsley2013convergence} that $\tau^{(\mu)} _{N_i}(\cB_{\frac 1r}(T_0))$ converges to $\tau^{(\mu)}((\cB_{\frac 1r}(T_0))$ as $i\to\infty$. Using Lemma \ref{lem:ballbound1}, this implies 
		\bb \nonumber
		\tau^{(\mu)}(\mathcal{B}_{\frac{1}{r}}(T_0) )\geq \Lambda(T_0,M)
		\ee
		for any $T_0\in\Omega_{\rm fin}$ and any $M >0$. Letting $M\to\infty$ we obtain 
		\bb \label{ineqeq}
		\tau^{(\mu)}(\mathcal{B}_{\frac{1}{r}}(T_0) )\geq \Lambda(T_0)\,.
		\ee
		Finally, using Lemma \ref{lem:addupto} and the fact that $\tau^{(\mu)}$ is a probability measure, it follows that equality holds in \eqref{ineqeq}. Since any Borel probability measure on $\cT$ is uniquely determined by its value on balls, by Theorem  2.2 in \cite{billingsley2013convergence}, this proves that the limit $\tau^{(\mu)}$ is unique. 
		
		That $\tau^{(\mu)}$ is supported on $\cT_\infty$ is clear by construction, since each finite tree is isolated in $\cT$. That $\tau^{(\mu)}$ is supported on $\Omega$ follows similarly from the fact that $\Omega$ is a closed subset of $\cT$ and hence its complement is a countable union of balls with vanishing measure. This shows that $\tau^{(\mu)}$ is supported on $\Omega_\infty$. 
	\end{proof}
	
	\subsection{Properties of the local limit}\label{sec:3.3}
	In this section we briefly provide a description of the measures $\tau^{(\mu)},\, \mu<\mu_0$, in terms of  branching processes  \cite{athreyaney1972branching, abraham2015introduction}. 
	
	Recall, that a Bienaymé-Galton-Watson (BGW) branching process (with one type of individual) is defined in terms of an offspring probability distribution $p(n), n=0,1,2,\dots,$ fulfilling 
	\bb\label{offspring1} \nonumber
	\sum_{n=0}^\infty p(n) =1\,.
	\ee
	The corresponding BGW tree is the probability measure $\lambda$ on $\cT$ defined by setting
	\bb\label{defGWtree1}
	\lambda(\cB_{\frac 1r}(T)) = \prod_{v\in \cup_{s=1}^{r-1}D_s(T)} p(\sigma(v)-1)\,,
	\ee
	for any $T\in\cT$. 
	It is well known that \eqref{defGWtree1} defines a Borel probability measure $\lambda$ on $\cT$ and that $\lambda$ is supported on $\cT_{\rm fin}$ if and only if  $p$ is subcritical or critical, i.e. if the average offspring $m$ satisfies 
	$$
	m :=\sum_{n=0}^\infty np(n) \leq 1\,.
	$$
	
	Given $p$ as above, the relevant branching process in our case has two types of individulas, called \emph{normal} and \emph{special}, whose offspring probabilities are restricted such that normal individuals can have $n=0, 1,2,\dots$ normal offspring with probability $p(n)$, but no special offspring, while special individuals can have $n=1,2,3,\dots$ off-spring, of which the first one is special, while the remaining $n-1$ are normal with probability $p(n-1)$. Moreover, the unique initial individual is special with a unique (special) decendant.
	By standard arguments, this defines a Borel probability meaure $\hat\lambda$ on $\Omega_\infty$  by setting
	\bb\label{defGWtree2}
	\hat\lambda(\cB_{\frac 1r}(T)) = \prod_{s=1}^{r-1}\Big(\prod_{v\in D_s(T)\setminus\{w_s\}} p(\sigma(v)-1)\Big)\cdot p(\sigma(w_s)-2)\,,
	\ee
	for any $T\in\Omega_\infty$, where $w_s$ denotes the leftmost vertex in $D_s(T)$. Obviously, the vertices corresponding to special individuals span the leftmost infinite path  $\omega_T$ starting from the root of $T$. If $p$ is subcritical or critical, the subtrees spanned by the descendants of any normal vertex is finite with probability $1$ and it can be shown that with $\hat\lambda$-probability $1$, the leftmost path is the unique infinite path in $T$ starting from the root, i.e. $T$ has a single spine equal to $\omega_T$, and $T$ is obtained by grafting finite trees at the spine vertices to the right of $\omega_T$. Due to the multiplicative structure of \eqref{defGWtree2}, these branches are independent with identical distribution equal to the subcritical or critical BGW tree with offspring probability $p$, while the degrees $\sigma(v)$ of the spine vertices $v$ (different from the root)  are likewise independently and identically distributed according to $p(\sigma(v)-2), \sigma(v)\geq 2$.
	In the following we shall denote single spine trees distributed according to $\hat\lambda$ by $\hat T$.  
	
	The interpretation of $\tau^{(\mu)}$ as a BGW tree is given by the following proposition.
	
	\begin{prop} \label{corconditional}
		For $\mu<\mu_0$, the measure $\tau^{(\mu)}$  equals the BGW measure $\hat\lambda$ defined above corresponding to the subcritical BGW tree with offspring probabilities given by
		\bb \label{branchingprob}
		p(n) = X(g_c(\mu))^n \big(1-X(g_c(\mu))\big) = e^{\mu n}(1-e^\mu)\,,~~ n=0,1,2,\dots .
		\ee
	\end{prop}
	\begin{proof}
		Since $X(g_c(\mu))<\frac 12$, it is clear that  $p$ given by \eqref{branchingprob} defines a probability distribution with mean 
		\bb\label{mean}
		m=\frac{X(g_c(\mu))}{1-X(g_c(\mu))} < 1\,,
		\ee
		and hence defines a subcritical BGW process. Using \eqref{branchingprob} in \eqref{defGWtree2}, it is seen for a given $T_0\in \Omega_{\rm fin}$ of height $r$ and with $|D_r(T_0)|=K$ that
		\begin{align}\nonumber
			\hat \lambda (\cB_{\frac 1 r} (T_0)) 	&= X(g_c(\mu))^{|E(T_0)|-r}\big(1-X(g_c(\mu))\big)^{|E(T_0)|-K}\\\nonumber& = g_c(\mu)^{|E(T_0)|-K}X(g_c(k))^{K-r}\\& = e^{-\mu(r-1)}\cdot g_c(k)^{|E(T_0)|-K}X(g_c(\mu))^{K-1}\,, \nonumber
		\end{align}
		where  \eqref{eq:XXm} has been used in the second step and \eqref{eq:gc} in the third step. The last expression is seen to coincide with \eqref{limit1}, which completes the proof of the proposition.	\end{proof}
	
	Letting $\mathbb E_\mu$ denote the expectation w.r.t. $\tau ^{(\mu)} $, the following result on the average growth of balls around the root of $\hat T$ as a function of radius is an easy consequence of familiar results about BGW processes.
	
	\begin{cor}\label{growth1} For $\mu<\mu_0$, the following asymptotic relations hold, 
		\begin{eqnarray}
			\mathbb E_{\mu}(|D_r|) &=&  \frac{1}{1-m} + O(m^{r-1})\,,\label{growthD1av}\\
			\mathbb E_{\mu}(|B_r|) &=&  \frac{r}{1-m}+ O(1)\label{growthB1av}\,,
		\end{eqnarray}
		where $m$ is given by \eqref{mean}.
	\end{cor}
	\begin{proof}  Let $T_s$  be the branches of $\hat T$ grafted at the spine vertex at distance $s$ from the root. It follows from Proposition \ref{corconditional} and the preceding remarks that these are i.i.d. according to the subcritical BGW tree with offspring distribution given by \eqref{branchingprob}. In particular, we have (see e.g. Ch. 1 of \cite{athreyaney1972branching}) 
		$$
		\mathbb E_\mu(|D_r(T_s)|) = m^{r-1}\,,\quad r\geq 1\,.
		$$
		Since 
		\bb\label{decompDr1} \nonumber
		|D_r(\hat T)| = 1 + \sum_{s=1}^{r-1} |D_{r-s+1}(T_s)|\,,
		\ee
		it follows  that 
		$$
		\mathbb E_\mu(|D_r(\hat T)|) = 1+ m\frac{1-m^{r-1}}{1-m}\,,
		$$
		from which the first relation follows. Using 
		\bb\label{DB}
		|B_r(\hat T)| = \sum_{s=1}^r |D_s(\hat T)|\,,
		\ee
		eq. \eqref{growthB1av} follows from \eqref{growthD1av}. 	\end{proof}
	
	A corresponding result on a.s. asymptotic growth of individual trees is obtained by arguments essentially identical to those given for the corresponding  result in \cite{meltem2021height}. 
	
	\begin{prop}\label{growth2}
		Let $\mu<\mu_0$. There exist  constants $C_1, C_2>0$ and for $\tau^{(\mu)}$-a.e.  $\hat T$ a number $r_0(\hat T) \in\mathbb N$, such that 
		\begin{align}\label{eq:growth2} \nonumber
			1&\leq |D_r(\hat T)| \leq C_1 \cdot \ln r\quad \mbox{and}\quad 
			r\leq |B_r(\hat T)| \leq C_2 \cdot r\ln r
		\end{align}
		for all $r\geq r_0(\hat T)$.
	\end{prop}

	\begin{remark}
		The volume growth exponent $d_h$ of a tree $T\in\cT_\infty$, defined by
		\bb\label{defHausdorff} \nonumber
		d_h := \lim_{r\to\infty} \frac{\ln |B_r(T)|}{\ln r}\,,
		\ee
		is commonly referred to as the Hausdorff dimension of $T$, provided the limit exists. Proposition \ref{growth2} shows that $d_h=1$ $\tau^{(\mu)}$-a.s. for $\mu<\mu_0$.
	\end{remark}

	\section{The case $\mu=\mu_0$}\label{sec:4}
	
	This case can be treated by techniques similar to those found in \cite{meltem2021alpha}, starting with the determination of the singular behaviour of the partition function as given by the following proposition, whose somewhat technical proof is provided in the appendix. In the following, we let for each $a>0$ the wedge $V_a$ be given by
	$$
	V_a= \{z\in\mathbb C\mid |{\rm Im}z|\,<\, a{\rm Re}z\}\,.
	$$

	\begin{prop}\label{prop:ZNln2}
		There exists a constant $c_0\in\mathbb R$ such that for any given $a>0$ we have 
		\bb\label{ZNln2}
		W^{(\mu_0)} (g) = -\frac 12 \ln(1-4g) + c_0 + O(|1-4g|^{\frac 15})
		\ee
		for all complex $g$ close to $\frac 14$ satisfying $\sqrt{1-4g}\in V_a$.
	\end{prop}
	
	Using transfer theorems \cite{flajolet2009analytic}, this implies the following asymptotic behaviour of $W^{(\mu_0)}_N$.
	
	\begin{cor} For $N$ large, it holds that
		$$
		W^{(\mu_0)}_N = \frac{1}{2N}\cdot 4^N\big(1+ o(1)\big)\,.
		$$
	\end{cor}
	Using  this result the arguments of section \ref{sec:3} can be repeated to establish the following theorem.
	
	\begin{thm} The measures $\tau^{(\mu_0)}_N$ converge weakly to the probability measure $\tau^{(\mu_0)}$ on $\Omega_\infty$ defined by
		$$
		\tau^{(\mu_0)}(\cB_{\frac 1r}(T_0)) = 4^{-|T_0|} 2^{K+r}\,,
		$$
		for any tree $T_0\in \Omega$ of height $r$ and with $K$ vertices at height $r$. 
		
		This measure equals the BGW measure corresponding to the BGW process with two types of individuals defined by the offspring probabilities $p(n) = 2^{-(n+1)}, n=0,1,2,\dots$, as in section \ref{sec:3.3}.  
	\end{thm}
	
	In particular, $\tau^{(\mu_0)}$ is supported on single-spine trees with finite branches grafted to the right that are i.i.d. according to the BGW measure $\rho$  defined in \eqref{defGWtree1} corresponding to the offspring distribution $p(n)$ as given. For later use we note that $\rho$ is equivalently given by
	\bb\label{def:rho}
	\rho(T) = 2\cdot 4^{-|T|}\,,\quad T\in \cT_{\rm fin}\,.
	\ee
	
	The proof of Corollary \ref{growth1} with $m=1$ now yields the following result.
	
	\begin{cor}\label{growth1c} The following asymptotic relations hold, 
		\begin{eqnarray} \nonumber
			\mathbb E_{\mu_0}(|D_r|) &=&  r\,,\label{growthDav}\\ \nonumber
			\mathbb E_{\mu_0}(|B_r|) &=&  \frac 12 r^2+ O(r)\label{growthBav}\,.
		\end{eqnarray}
	\end{cor}
	
	Moreover, by a close analogue to the arguments applied for the UIPT (see e.g. Proposition 1 of \cite{durhuus2010spectral}) one obtains  
	
	\begin{prop}\label{growth2c}
		There exist  constants $C_1, C_2>0$ and for $\tau^{(\mu_0)}$-a.e.  $\hat T$ a number $r_0(\hat T) \in\mathbb N$, such that 
		\begin{align}\label{eq:growth2c} \nonumber
			C_1\cdot (\ln r)^{-2}r^2\; \leq \; |B_r(\hat T)|\; \leq \;C_2 \cdot (\ln r)r^2
		\end{align}
		for all $r\geq r_0(\hat T)$.
	\end{prop}

	\section{The case $\mu>\mu_0$}\label{sec:5}
	
	\subsection{The partition function}\label{sec:5.1}
	
	In this subsection we determine the asymptotic behaviour of the partition functions $W_N^{(\mu)}$ for large $N$, in case $\mu>\mu_0$. 
	
	\begin{thm}\label{thm:4.1} For each  $\mu> \mu_0$ it holds for any $\delta\in ]0,\frac 16[$ that
		\bb\label{eq:asympZ}
		W^{(\mu)}_N = 4 \cdot e^\mu  \sqrt{\frac{\pi}{B}}\frac{\mu+\ln 2}{2} e^{-A N^{\frac 13}}N^{-\frac 56}4^N\Big(1+O(N^{-\delta})\Big) \,,
		\ee
		for $N$ large, where
		$$
		A =  3\Big(\frac{\pi(\mu+\ln 2)}{2}\Big)^{\frac 23} \quad\mbox{and}\quad B = 3\Big(\frac{(\mu+\ln 2)^2}{4\pi}\Big)^{\frac 23} \,.
		$$
	\end{thm}
	
	\begin{proof} In \cite{meltem2021height} we established a similar formula for the asymptotic behaviour of the partition function $Z^{(\mu)}_N$ given by \eqref{eq:partZ}
		by using a rather standard saddlepoint approximation. It was shown that the contribution of the terms in \eqref{AmN} corresponding to $k\geq 2$ is suppressed by a factor of the form $e^{-cN^{\frac 13}}$ relative to the contribution from the $k=1$ term. Rewriting the contribution from the $k=1$ terms as 
		\bb\label{deftildeZN}
		\tilde Z^{(\mu)}_N := 4^N\sum_{m=2}^N \frac{e^\mu}{m+1}\tan^2\frac{\pi}{m+1}\; e^{-f_N(m+1)}\,,
		\ee
		where 
		\bb\label{eq:deffN}  \nonumber
		f_N(t) = \mu t + N\ln\Big( 1+\tan^2\frac{\pi}{t}\Big)\,,\quad t>2\,,
		\ee
		the saddlepoint $t_0$ is determined by $f_N'(t_0)=0$ and is given by 
		\bb\label{eq:minasymp} \nonumber
		t_0 = \Big(\frac{2\pi^2 N}{\mu}\Big)^{\frac 13} + O(N^{-\frac 13})\,.
		\ee
		Moreover, the dominant contribution to the sum \eqref{deftildeZN} originates from values of $m$ in $[t_0-N^{\frac 13-\delta^\prime}, t_0+N^{\frac 13-\delta^\prime}]$, where $0<\delta^\prime<\frac 16$. Since in this interval we have
		$$
		\Big( 1+\tan^2\frac{\pi}{m+1}\Big)^{\frac{m-1}{2}} = 1+O(N^{-\frac 13})\,,
		$$
		a comparison of the expressions  \eqref{BmN} and \eqref{AmN} shows that the saddlepoint contribution to $W^{(\mu)}_N$ is identical to the one for $2\tilde Z_N ^{(\mu+\ln 2)}$ up to a factor $1+O(N^{-\frac 13})$, which is in fact identical to the right-hand side of \eqref{eq:asympZ} with $\delta=3\delta^\prime-\frac 13$. We leave it to the reader to check that the remaining contributions to $W^{(\mu)}$ are subleading, provided $\delta^\prime>\frac 19$, by arguments essentially identical to the ones in \cite{meltem2021height}. This proves that \eqref{eq:asympZ} holds.
	\end{proof}
	
	\subsection{Ball volumes and the local limit}\label{sec:5.2}
	We next calculate the asymptotic behaviour of the $\tau^{(\mu)}_N$-volume of balls that will allow us to prove tightness of the sequence $(\tau^{(\mu)}_N)$, and also to show weak convergence.
	
	\begin{lem}\label{lem:4.2} Assume $\mu>\mu_0$ and that  $T_0\in \Omega^{(r)}$ and set $K = |D_r(T_0)|$.  
		Given  $0<\epsilon<\frac 1K$ and $M\in\mathbb N$, it holds for any $\delta\in ]0,\frac 16[$ that 
		\begin{align}
			\tau^{(\mu)}_N(\cB_{\frac 1r}(T_0)) ~ \geq ~ &\frac{e^{-\mu(r-1)}}{4^{|T_0|-K}} \sum_{R=1}^{K} \binom{K-1}{R-1}\frac{1}{(R-1)!}\Big(\frac{\mu + \ln 2}{2}\Big)^{R-1}\\ & \Big(\sum_{S=1}^M C_{S-1}4^{-S}\Big)^{K-R}(1-\epsilon K)^K(1+O(N^{-\delta})) \label{est1}
		\end{align}
		for $N$ large.  
	\end{lem}
	\begin{proof}
		As noted earlier, the elements $T$ of $\cB_{\frac 1r}(T_0)$ are obtained by grafting $K$ trees $T_1\dots, T_K$ onto $T_0$ at the $K$ vertices of maximal height, see Fig. 2. We call $T_1,\dots, T_K$ the branches of $T$. For $T\in\Omega_N$ we then have 
		$$
		N= \sum_{i=1}^K |T_i| + |T_0| - K\,,
		$$
		and we shall call the branch $T_i$ \emph{small} if $|T_i|\leq M$, while we call it \emph{large} if $|T_i| > \epsilon N$.  Note that, if $N>\frac{M}{\epsilon}$, no $T_i$ can be both small and large, and if also $N>K\cdot M +|T_0|-K$ there must be at least one large branch $T_i$. Moreover, if $h(T)> M+r-1$ then $h(T_1)>M$, which in particular ensures that $T_1$ is not small. Assuming $T$ fulfills these conditions in the following, let $\Omega_{T_0,R}\,, 1\leq R\leq K$, denote the subset of $\cB_{\frac 1r}(T_0)\cap \omega_N$ consisting of trees $T$ whose large branches are precisely $T_1,\dots, T_R$. Since $\tau^{(\mu)}_N$ restricted to $\cB_{\frac 1r}(T_0)$ is invariant under permutation of the branches different from $T_1$, we have
		\bb\label{Rsum}
		\tau^{(\mu)}_N(\cB_{\frac 1r}(T_0)) = \sum_{R=1}^{K} \binom{K-1}{R-1} \tau^{(\mu)}_N(\Omega_{T_0,R})\,.
		\ee
		Hence, we proceed to estimate   $\tau^{(\mu)}_N(\Omega_{T_0,R})$. 
		
		Denoting $|T_i|$ by $N_i$, we have
		\begin{align}\nonumber
			\tau^{(\mu)}_N(\Omega_{T_0,R})\,  \geq\,  & \frac{e^{-\mu(r-1)}}{ W^{(\mu)}_N}  \sum_{m=M+1}^\infty \sum_{N_{R+1},\dots, N_K\leq M} C_{N_{R+1}-1}\dots C_{N_K-1}\\ \nonumber
			& ~ \sum_{\substack{N_1+\dots +N_R \\ = N-N_{R+1}-\dots -N_K+K-|T_0|\\ N_1,\dots, N_R > \epsilon N}} e^{-\mu m} B_{m,N_1}A_{m,N_2}\dots A_{m,N_{R}}\, .
		\end{align}  
		Now, consider the last sum for fixed values of $N_{R+1},\dots ,N_K$ and set 
		$$
		N' = N_{R+1}+\dots+N_K +|T_0|-K\,.
		$$
		A lower bound on this sum is obtained by restricting the sum to values of $m$ fulfilling $|m-t_0|\leq N^{\frac 13-\delta^\prime}$, where $\delta^\prime > 0$.  Hence, setting 
		$$ x_m = \frac{m}{N^{\frac 13}}\quad \mbox{and}\quad x_0= \frac{t_0}{N^{\frac 13}}\,,$$
		we proceed to further estimate this lower bound
		\begin{align}\label{defVN}
			&V_N := \sum_{|x_m-x_0|\leq N^{-\delta^\prime}} \sum_{\substack{N_1+\dots +N_R =N-N'\\ N_1,\dots, N_R > \epsilon N}} e^{-\mu m}B_{m,N_1}\prod_{s=2}^R A_{m,N_s}\,.
		\end{align} 
		
		Recalling the expression \eqref{AmN} for $A_{m,N}$, we write
		\bb\label{eq:A1}
		A_{m,N_s} =  4^{N_s}\frac{1}{m+1}\tan^2\frac{\pi}{m+1}\big(1+ \tan^2\frac{\pi}{m+1}\big)^{-N_s} \Big(1+\sum_{k=2}^{\lfloor\frac m2\rfloor}\frac{\tan^2\frac{\pi k}{m+1}}{\tan^2\frac{\pi}{m+1}}\Big(\frac{1+\tan^2\frac{\pi}{m+1}}{1+\tan^2\frac{\pi k}{m+1}}\Big)^{N_s}\Big)\,.
		\ee
		Using that  
		$$
		\frac{1+\tan^2\frac{\pi}{m+1}}{1+\tan^2\frac{\pi k}{m+1}} \leq \frac{1+\tan^2\frac{\pi}{m+1}}{1+\tan^2\frac{2\pi }{m+1}} = e^{-\frac{3\pi^2}{(m+1)^2}\big(1+O(\frac{1}{m^2})\big)}
		$$
		for $2\leq k\leq \lfloor\frac m2\rfloor$, the sum in \eqref{eq:A1} is bounded from above by
		$$
		(m+1)^3 e^{-\frac{3\pi^2(N_s-1)}{(m+1)^2}\big(1+O(\frac{1}{m^2})\big)}\,. 
		$$
		In the range of $m$ considered and for $N_s > \epsilon N$, we hence get that 
		\bb\label{eq:A2} \nonumber
		A_{m,N_s} =  4^{N_s}\frac{1}{m+1}\tan^2\frac{\pi}{m+1}\big(1+ \tan^2\frac{\pi}{m+1}\big)^{-N_s}\big(1+ O(e^{-\epsilon c N^{\frac 13}})\big)\,,
		\ee
		where $c>0$ is a numerical constant. A similar estimate obviously holds for $B_{m,N}$ (up the the factor $2^{-(m-1)}$), such that the product in \eqref{defVN} takes the form 
		\begin{align}\label{eq:prod1}
			B_{m,N_1}\prod_{s=2}^R A_{m,N_s} &=  4^{N-N'}2^{-(m-1)}\Big(\frac{1}{m+1}\tan^2\frac{\pi}{m+1}\Big)^R\big(1+ \tan^2\frac{\pi}{m+1}\big)^{-(N-N')}\nonumber\\ &~~~~~~~~~~~~~~~~~~~~~~~~~~~~~~~~~~~~~~~~~~~~\cdot \big(1+ O(e^{-\epsilon c N^{\frac 13}})\big)^R\,.
		\end{align}
		Inserting \eqref{eq:prod1} into \eqref{defVN}, the sum over $N_1,\dots, N_R$ can be performed and yields a combinatorial factor 
		\begin{align}\nonumber
			\binom{N-N'-R\lfloor \epsilon N\rfloor +R-1}{R-1} & \geq \frac{N^{R-1}}{(R-1)!}\Big(1-\frac{N'+R\lfloor\epsilon N\rfloor}{N}\Big)^{R-1} \\ &\geq  \frac{N^{R-1}}{(R-1)!}(1- \epsilon R)^{R-1}\big(1+O(N^{-1})\big)\,.\nonumber
		\end{align}
		Thus, we obtain
		\begin{align}\nonumber
			V_N \geq 2 \cdot 4^{N-N'}\frac{1}{(R-1)!}\sum_{|x_m-x_0|\leq N^{-\delta}}& N^{R-1} e^{-(\mu+\ln 2) m}\Big(\frac{1}{m+1}\tan^2\frac{\pi}{m+1}\Big)^R \\ & \cdot\big(1+ \tan^2\frac{\pi}{m+1}\big)^{-(N-N')} (1-\epsilon K)^{K} \big(1+O(N^{-1})\big)\,.\nonumber
		\end{align} 
		The sum can now be estimated by repeating the arguments leading to  \eqref{eq:asympZ}, 
		and we arrive at 
		\bb\nonumber
		V_N \geq 4e^\mu \cdot 4^{N-N'}\frac{1}{(R-1)!}\sqrt{\frac{\pi}{B}}\Big(\frac{\mu+\ln 2}{2}\Big)^R e^{-A N^{\frac 13}}N^{-\frac 56} \cdot (1-\epsilon K)^{K} \big(1+O(N^{\frac 13-3\delta^\prime})\big)\,,
		\ee
		provided $\delta^\prime >\frac 19$. 
		Using this estimate as well as Proposition \ref{thm:4.1}, we obtain
		\begin{align}\nonumber
			&\tau^{(\mu)}_N(\Omega_{T_0,R})\,  \geq\,  e^{-\mu(r-1)}\frac{V_{N}}{ W^{(\mu)}_N}\,\geq\, \frac{e^{-\mu(r-1)}}{4^{|T_0|-K}} \frac{1}{(R-1)!}\Big(\frac{\mu+\ln 2}{2}\Big)^{R-1} \Big(\sum_{S=1}^M C_{S-1}4^{-S}\Big)^{K-R}\\ \nonumber &~~~~~~~~~~~~~~~~~~~~~~~~~~~~~~~~~~~~~~~~~~~~~~~~~~~~~~~~~~~~~\cdot (1-\epsilon K)^K(1+O(N^{\frac 13 -3\delta^\prime}))\,,
		\end{align}
		from which \eqref{eq:asympZ} follows by using \eqref{Rsum} and choosing $\frac 19<\delta^\prime<\frac 16$.
	\end{proof}
	
	In the following, let $\Xi(T_0;M,\epsilon)$ denote the large-$N$ limit of the right-hand side of \eqref{est1},
	$$
	\Xi(T_0;M,\epsilon) =  \frac{e^{-\mu(r-1)}}{4^{|T_0|-K}} \sum_{R=1}^{K} \binom{K-1}{R-1}\frac{1}{(R-1)!}\Big(\frac{\mu+\ln 2}{2}\Big)^{R-1}\Big(\sum_{S=1}^M C_{S-1}4^{-S}\Big)^{K-R}(1-\epsilon K)^K\,,
	$$
	and set 
	\bb\label{defXi}
	\Xi(T_0) := \lim_{\substack{\epsilon\to 0\\ M\to \infty}} \Xi(T_0,M,\epsilon) =  \frac{e^{-\mu(r-1)}}{4^{|T_0|-K}} \sum_{R=1}^{K} \binom{K-1}{R-1}\frac{1}{(R-1)!}\Big(\frac{\mu+\ln 2}{2}\Big)^{R-1}2^{R-K}\,,
	\ee
	where we have also used that 
	$$
	\sum_{S=1}^\infty C_{S-1} 4^{-S} = X(\frac 14) = \frac 12\,.
	$$
	We then have following analogue of Lemma \ref{lem:addupto}.
	
	\begin{lem}\label{lem:sumrule2} For all $r\geq 1$ it holds that
		\bb\label{sumrule2}
		\sum_{T_0\in\Omega^{(r)}} \Xi(T_0) = 1\,.
		\ee
	\end{lem}
	
	\begin{proof} 
		We use induction with respect to $r$. The case $r=1$ being trivial, let $r\geq 2$ and assume \eqref{sumrule2} holds for $r-1$. 
		It is convenient to rewrite \eqref{defXi} as 
		$$
		\Xi(T_0) = \sum_{\substack{B\subseteq D_r(T_0)\\ w_r\in B}} \Xi(T_0,B)\,,
		$$
		where $w_r$ is the leftmost vertex at graph distance $r$ from the root and 
		$$
		\Xi(T_0,B) =  \frac{e^{-\mu(r-1)}}{4^{|T_0|-|D_r(T_0)|}}\frac{1}{(|B|-1)!}\Big(\frac{\mu+\ln 2}{2}\Big)^{|B|-1}2^{-|D_r(T_0)\setminus B|}\,,
		$$
		such that the sum in \eqref{sumrule2} becomes a double sum over $(T_0,B)$. 
		
		Let  $T'_0\in\Omega^{(r-1)}$ be a fixed tree with $K'$ vertices at height $r-1$. Given a subset $B'\subseteq D_{r-1}(T_0')$ containing $w_{r-1}$, we then consider the contribution to the sum in \eqref{sumrule2} from all $(T_0,B)$ such that $T_0\in \Omega^{(r)}$ coincides with $T_0'$ up to height $r-1$ and such that $B'$ is the set of parents to vertices in $B$, i.e. $\phi_r(B)=B'$ in the notation of section \ref{sec:2}, where $\phi_r$ is the $r$'th parent map of $T_0$. We claim that this contribution is precisely  $\Xi(T'_0,B')$.
		
		In order to establish the claim, we first note that for any ordered set $B$ as above with $|B|=R$ the number of surjective order preserving maps $\varphi: B\to B'$ equals $ \binom{R-1}{ R'-1}$, where $R'=|B'|$, since necessarily $\phi_r(w_r)=w_{r-1}$. Moreover, for any given such $\varphi$, the number of ways of extending it to an order preserving map $\phi_r$ (not necessarily surjective) from any ordered set of $K$ elements, containing $B$ as an ordered subset, into $D_{r-1}(T_0')$ such that the first element is mapped to $w_{r-1}$ is easily seen to equal $\binom{ K+K'-2}{ K'+R-2}$. Hence, we have 
		\begin{align}\nonumber
			\sum_{\substack{(T_0,B): T_0\in\Omega^{(r)},\\ w_r\in B\subseteq D_r(T_0)\\
					B_{r-1}(T_0)=T_0',\, \phi_r(B)=B'}} \Xi (T_0,B) &=  \frac{e^{-\mu(r-1)}}{4^{|T_0'|}} \Big[ \sum_{K=1}^\infty\sum_{R=1}^K  \binom{K+K'-2}{K'+R-2} \binom{R-1}{R'-1} \\ \nonumber &~~~~~~~~~~~~~~~~ \Big(\frac{\mu+\ln 2}{2}\Big)^{R-1}\frac{1}{(R-1)!}\,2^{R-K}~\Big] \,. 
		\end{align}
		Interchanging the summation order and using the identity \eqref{combid},
		the right-hand side becomes 
		\begin{align}\nonumber
			\frac{e^{-\mu(r-1)}}{4^{|T_0'|}} \sum_{R=R'}^\infty \binom{R-1}{R'-1}\frac{(\mu+\ln 2)^{R-1}}{(R-1)!}\,2^{K'}
			&=  \frac{e^{-\mu(r-1)}}{4^{|T_0'|}} \sum_{R=R'}^\infty \frac{(\mu+\ln 2)^{R-1}}{(R'-1)!(R-R')!}2^{K'}\\\nonumber
			&=  \frac{e^{-\mu(r-2)}}{4^{|T_0'|}} \frac{(\mu+\ln 2)^{R'-1}}{(R'-1)!}\,2^{K'+1}\,,
		\end{align}
		which is seen to be equal to $\Xi(T'_0,B')$, as claimed. Hence we have
		
		\begin{align}\nonumber
			\sum_{\substack{(T_0,B): T_0\in\Omega^{(r)}\\ w_r\in B\subseteq D_r(T_0)}} \Xi(T_0,B) &= \sum_{\substack{(T'_0,B'):T'_0\in\Omega^{(r-1)}\\ w_{r-1}\in B^\prime\subseteq D_{r-1}(T'_0)}} \sum_{\substack{(T_0,B): T_0\in\Omega^{(r)}\\ w_r\in B\subseteq D_r(T_0)\\
					B_{r-1}(T_0)=T_0',\, \phi_r(B)=B'}} \Xi(T_0,B)\\& =  \sum_{\substack{(T'_0,B'): T'_0\in\Omega^{(r-1)}\\ w_{r-1}\in B'\subseteq D_{r-1}(T'_0)}} \Xi(T'_0,B') =1\,,\nonumber
		\end{align}
		where the induction assumption has been used in the last step. This finishes the proof.
	\end{proof}
	
	The main result on the existence of the local limit can now be established. 
	
	\begin{thm}\label{thm:limit2}
		For $\mu>\mu_0$, the sequence $(\tau^{(\mu)}_N)$ is weakly convergent to a Borel probability measure $\tau^{(\mu)}$ on $\cT$ characterized by
		\bb\label{ballmeasure2}
		\tau^{(\mu)}(\cB_{\frac 1r}(T_0)) = \Xi(T_0) = \frac{e^{-\mu(r-1)}}{4^{|T_0|}}\,2^{K+1}\,\sum_{R=1}^{K} \binom{K-1}{R-1}\frac{(\mu+\ln 2)^{R-1}}{(R-1)!}\,,
		\ee 
		for any tree $T_0\in\Omega^{(r)}$, where $K=|D_r(T_0)|$ and $r\geq 1$. Moreover, the measure is concentrated on $\Omega_\infty$ and in particular $\tau^{(\mu)}(\cB_{\frac{1}{r}}(T_0))$ vanishes unless $T_0$ is one-sided.
	\end{thm}
	
	\begin{proof} This follows by the same line of reasoning as in the proof of Theorem \ref{thm:limit1}, using Lemmas \ref{lem:sumrule2} and \ref{lem:4.2}. We leave the details to the reader. 
	\end{proof}
	
	\begin{remark}\label{numuchar} The formula corresponding to \eqref{ballmeasure2} characterising the measure $\nu^{(\mu)}\,, \mu>0$, states that, for arbitrary $T_0\in\cT^r$ and with $K=|D_r(T_0)|$, we have
		\bb\label{ballmeasure3}
		\nu^{(\mu)}(\cB_{\frac 1r}(T_0)) = \frac{e^{-\mu(r-1)}}{4^{|T_0|}}\,2^{K+1}\,\sum_{R=1}^{K} \binom{K}{R}\frac{\mu^{R-1}}{(R-1)!}\,,
		\ee 
		see Theorem 4.5 in \cite{meltem2021height}. This result is needed in the subsequent proof.
	\end{remark}
	
	\smallskip
	
	Given $n \in \mathbb N$, let us define the standard $(n-1)$-simplex by
	
	$$
	\Delta_{n} := \{ (x_1,\dots,x_n) \mid x_1 + \dots +x_n = 1\,,\; x_1,\dots, x_n > 0 \}\,,
	$$
	as well as the scaled simplex
	
	$$
	\mu\cdot\Delta_{n} := \{ (\mu_1,\dots,\mu_{n}) \mid \mu_1 + \dots +\mu_{n} = \mu, ~\mu_1, \dots \mu_{n} >0 \} \,,\;\mbox{for $\mu>0$}.
	$$
	
	In the next subsection we shall need the following decomposition result, where $d\omega_n$ denotes the normalized Lebesgue measure on the simplex $(\mu+\ln 2)\cdot\Delta_n$ and the reader should recall that the BGW measure $\rho$ is defined by \eqref{def:rho}.
	
	\begin{prop} \label{propdecomp}
		Assume $\mu>- \ln 2$ and that $T_0 \in \Omega^{(r)}$ and set $K=|D_r(T_0)|$. Under the identification $\cB_{\frac 1r}(T_0) = \cT^{D_r(T_0)}$ via grafting trees onto $T_0$ at the vertices in $D_r(T_0)$ (see Remark \ref{remarkgraft}), it holds for arbitrary trees $T_1, \dots , T_K \in \cT_{\text{fin}}$, with $h_i=h(T_i)$, that 		
		\begin{align} \label{decomposition} \nonumber 
			\tau ^{(\mu)} \mid_{\cB_{\frac 1r} (T_0)} (\cB_{\frac{1}{h_1}}(T_1) \times \dots \times \cB_{\frac{1}{h_K}}(T_K)) &= \sum\limits_{\substack{D \subseteq D_r(T_0)\\ w_r \in D}} \frac{e^{-\mu(r-1)}}{4^{|T_0|}} 2^{K+1} \frac{(\mu + \ln  2)^{|D|-1}}{(|D|-1)!} \\ \cdot\Big[ \int\limits_{(\mu+ \ln 2 )\cdot\Delta_{|D|}} d\omega_{|D|}\,   \tau^{(\mu_1 - \ln 2)} (\cB_{\frac{1}{h_1}} (T_1))    & \prod\limits_{i \in D\setminus\{w_r\}}  \nu^{(\mu_i)}(\cB_{\frac{1}{h_i}}(T_i)) \prod\limits_{j \notin D }\rho (\cB_{\frac{1}{h_j}} (T_j)) \Big]\,,
		\end{align}
		where by abuse of notation the $i$th vertex in $D_r(T_0)$ from the left, onto which $T_i$ is grafted, has been identified with the index $i$ for  $i= 1, \dots , K$; in particular, $1$ is used to denote the leftmost vertex $w_r$ in $D_r(T_0)$.
	\end{prop}
	
	\begin{proof}
		Setting $h= \max\{h_1, \dots , h_k\}$, the product set on the left-hand side of \eqref{decomposition} can be decomposed in a disjoint union of balls as 	
		$$
		\cB_{\frac{1}{h_1}}(T_1) \times \dots \times \cB_{\frac{1}{h_K}}(T_K) = \bigcup_{\substack{T \in \cT^{r+h-1} \\ gr(T_0;T_1,\dots,T_K) \subseteq T}} \cB_{\frac{1}{h+r-1}} (T) \,,
		$$
		\noindent where $gr(T_0;T_1,\dots,T_K) \subseteq T$ indicates that $T_0$ with $T_1,\dots, T_K$ grafted on its vertices at the maximal height is a subtree of $T$ with the same root edge. Thus, the left-hand side of \eqref{decomposition} can be written as a sum of the corresponding ball volumes. Using \eqref{ballmeasure2} and \eqref{ballmeasure3} together with the identity
		$$
		\int_{\Delta_n} d\omega \prod_{i=1}^nx_i^{m_i}  = \frac{(n-1)!}{(m_1+\dots+m_n +n-1)!}\prod_{i=1}^n m_i!\,,
		$$
		valid for non-negative integers $m_1,\dots, m_n$, the desired result follows in a straightforward way by suitably rearranging the mentioned sum. 
	\end{proof}
	
	\begin{remark}\label{spinesum} 	Noting that for any given $T\in \cT^{(m)}$ we have
		$$
		\cB_{\frac 1m}(T)\setminus\{T\} = \bigcup_{T^\prime\in\cT^{m+1}\cap\cB_{\frac{1}{m}}(T)}\cB_{\frac{1}{m+1}}(T^\prime)
		$$
		where the right-hand side is a disjoint union, it follows, using countable additivity on both sides that \eqref{decomposition} also holds if $\cB_{\frac{1}{h_i}}(T_i)$ is replaced by $\{T_i\}$ for one or more values of $i$. Since $\tau^{(\mu_1-\ln 2)}$ and $\nu^{(\mu)}$ are supported on infinite trees for $\mu_1, \mu>0$, it follows that the terms on the right-hand side vanish for those $D$ that contain one or more of those indices. In turn it likewise follows that \eqref{decomposition} also holds if $\cB_{\frac{1}{h_i}}(T_i)$ is replaced by $\cB^\infty_{\frac{1}{h_i}}(T_i) := \cB_{\frac{1}{h_i}}(T_i)\cap \cT_\infty$. It thus follows that the term on the right-hand side corresponding to a given $D$ equals the $\tau^{(\mu)}$-measure of the set of trees whose infinite branches are precisely the ones grafted on vertices in $D$. 
	\end{remark}
	
	\subsection{Properties of the local limit}\label{sec:5.3}
	
	In \cite{meltem2021height} we defined the probability measure $\tilde{\nu}^{(\mu)}$ on the set $\cT_s$ of spine trees by
	\bb \label{nutilde}
	\tilde{\nu} ^{(\mu)} (\cB^s_{\frac{1}{r}}(T^s))=e^{-(r-1)\mu} \frac{\mu^{R-1}}{(R-1)!} ~,~~ r \geq 1\,,
	\ee
	for any $T^s\in \cT_s$ with $R$ vertices at height $r$, and where 
	$$
	\cB^s_{a}(T^s) := \cB_a(T^s)\cap \cT_s
	$$ 
	denotes the ball in $\cT_s$ of radius $a>0$ around a spine tree $T^s$. Note that $\cT_s\subset \Omega$, and clearly any $T\in\Omega_\infty$ can be obtained in a unique way by grafting finite trees (branches) in any sector $(i,n)$ of $\chi(T)$ except in the first sector of the leftmost vertices at each height. In this way, we have by Remark \ref{remarkgraft} a homeomorphism (see Fig. \ref{figure:treespine})
	\bb\label{identification}
	G: \chi^{-1}(\cB^s_{\frac{1}{r}}(T^s)) \to \Omega_\infty \times \cT_\infty^{R-1}\times  \cT_{\rm fin}^\sigma\,,
	\ee
	where $R=|D_r(T^s)|$ is the number of leaves in $T_0^s := B_r(T^s)$ and 
	\bb\label{sigmaform}
	\sigma = 2|T^s_0| -R-r
	\ee
	is the total number of sectors associated with $T_0 ^s$, excluding those adjacent to the leaves and the first sector of the leftmost vertices.
	
	\begin{prop}\label{pushforward} Assume $\mu>\mu_0$ and let $\chi^*\tau^{(\mu)}$ denote the pushforward of $\tau^{(\mu)}$ under the spine map $\chi$. Then the following hold.
		
		i)\; $ \chi^*\tau^{(\mu)} = \tilde{\nu}^{(\mu+\ln 2)}$
		
		ii)\; The branches $T_{(i,n)}$ are i.i.d. according to the BGW-measure $\rho$.
	\end{prop}
	
	\begin{proof}
		As above, let $T^s$ be any spine tree with $|D_r(T^s)|=R$ and $B_r(T^s)= T^s _0$. Using the identification \eqref{identification}, we first evaluate
		\bb \label{moreballs}
		\tau^{(\mu)} \mid_{\chi^{-1}(\cB^s _{\frac{1}{r}}(T^s))} \Big( \cB^\infty_{\frac{1}{h_1}}(T_1) \times \dots \times \cB^\infty_{\frac{1}{h_R}}(T_R) \times \prod_{(i,n)} \{T_{(i,n)}\}\Big)
		\ee
		for $T_j,~T_{(i,n)}\in\cT_{\rm fin}$ fixed where $j=1,\dots,R$, $(i,n)$ labels the sector in which $T_{(i,n)}$ is grafted and $h_j=h(T_j)$.
		We further set $T^\prime_{(i,n)}=B_{r-|i|+1}(T_{(i,n)})$, where $|i|$ is the height of the vertex $i$ in $T^s_0$ and let $T_0$ be the tree of height $r$ obtained by grafting $T^\prime_{(i,n)}$ onto $T^s_0$ at $(i,n)$ instead of $T_{i,n}$ in each sector. If $l_{(i,n)} := |D_{r-|i|+1}(T_{(i,n)})|$, the number of vertices in $T_0$ at height $r$ equals $K= R+ \sum\limits_{(i,n)} l_{(i,n)}$. Together with \eqref{sigmaform} this implies
		$$
		4^{-|T_0|} 2^{K+1} = 2^{-(r-1)} \prod\limits_{(i,n)} 4^{-|T^\prime _{(i,n)}|} 2^{l_{(i,n)}+1} \, .
		$$
		
		\noindent By Proposition \ref{propdecomp} and Remark \ref{spinesum} this implies that \eqref{moreballs} equals
		\begin{align*}
			e^{-(\mu+ \ln2)(r-1)} \frac{(\mu + \ln 2)^{R-1}}{(R-1)!} & \Big[ \int\limits_{(\mu+\ln 2)\cdot\Delta_R}d\omega \cdot \tau^{(\mu_1 - \ln 2)} (\cB_{\frac{1}{h_1}} (T_1)) \prod\limits_{j=2}^R  \nu^{(\mu_j)}(\cB_{\frac{1}{h_j}}(T_j))\Big]  \\ &  ~~~~  \prod\limits_{(i,n)} 4^{-|T^\prime _{(i,n)}|} 2^{l_{(i,n)}+1}\prod_{k=1}^{l_{(i,n)}}\rho (\{T^k_{(i,n)}\})\,,
		\end{align*}
		where $T^k_{(i,n)},\, k=1,\dots,l_{(i,n)}$, are the branches grafted on the vertices at maximal height of $T^\prime_{(i,n)}$ to yield $T_{(i,n)}$ (provided $l_{(i,n)}\neq 0$).
		Upon noting that 
		$$
		4^{-|T^\prime _{(i,n)}|} 2^{l_{(i,n)}+1}\prod_{k=1}^{l_{(i,n)}}\rho (\{T^k_{(i,n)}\}) = \rho(T_{(i,n)})\,,
		$$
		we hence obtain 
		\begin{align*}
			&\tau^{(\mu)} \mid_{\chi^{-1}(\cB^s _{\frac{1}{r}}(T^s))} \Big( \cB^\infty_{\frac{1}{h_1}}(T_1) \times \dots \times \cB^\infty_{\frac{1}{h_R}}(T_R) \times \prod_{(i,n)} \{T_{(i,n)}\}\Big) \\ &= e^{-(\mu+ \ln2)(r-1)} \frac{(\mu + \ln 2)^{R-1}}{(R-1)!}\Big[ \int\limits_{(\mu+\ln 2)\cdot\Delta_R}d\omega \cdot \tau^{(\mu_1 - \ln 2)} (\cB_{\frac{1}{h_1}} (T_1)) \prod\limits_{j=2}^R  \nu^{(\mu_j)}(\cB_{\frac{1}{h_j}}(T_j))\Big] \prod\limits_{(i,n)}\rho (\{T_{(i,n)}\})\,.
		\end{align*}
		Summing the right-hand side over all finite trees $T_1,\dots, T_R$ of given \textit{fixed} heights $h_1,\dots,h_R$, respectively, as well as over all finite trees $T_{(i,n)}$  gives
		$$
		e^{-(\mu+ \ln2)(r-1)} \frac{(\mu + \ln 2)^{R-1}}{(R-1)!} = \tilde{\nu} ^{(\mu+\ln 2)} (\cB^s _{\frac{1}{r}} (T_0 ^s))\,,
		$$
		\noindent proving of the first statement of the proposition as well as 
		\begin{align*}
			&\tau^{(\mu)}\Big( \cB^\infty_{\frac{1}{h_1}}(T_1) \times \dots \times \cB^\infty_{\frac{1}{h_R}}(T_R) \times \prod_{(i,n)} \{T_{(i,n)}\}\,\Big|\,{\chi^{-1}(\cB^s _{\frac{1}{r}}(T^s))} \Big) \\ &= \Big[ \int\limits_{(\mu+\ln 2)\cdot\Delta_R}d\omega \cdot \tau^{(\mu_1 - \ln 2)} (\cB_{\frac{1}{h_1}} (T_1)) \prod\limits_{j=2}^R  \nu^{(\mu_j)}(\cB_{\frac{1}{h_j}}(T_j))\Big] \prod\limits_{(i,n)}\rho (\{T_{(i,n)}\})\,,
		\end{align*}
		which implies the second statement.
	\end{proof}
	
	\begin{figure}[H]
		\centering
		\begin{align*}
			\raisebox{0cm}{\treespine} 
		\end{align*}
		
		\captionof{figure}{Structure of the elements in $\chi^{-1}(\cB_{\frac 13} ^s (T_0 ^s)) \cap \Omega$ where $T_0 ^s$  of height $h(T_0 ^s)=3$ with root vertex $i_0$ is shown in black. The green blob and the red blobs indicate the infinite branches from $\Omega_\infty$ and $\cT_\infty$ respectively. Blue blobs are independent and distributed according to the critical BGW measure $\rho$.}
		\label{figure:treespine}
	\end{figure}
	
	Based on this proposition one easily obtains by arguments similar to those of \cite{meltem2021height} the following results on the volume growth of balls w.r.t. $\tau^{(\mu)}$.
	
	\begin{cor} \label{expected}
		The following statements hold for $\mu>\mu_0$.
		\begin{itemize}
			\item [i)] $\mathbb E_{\mu}(|D_r|)= (\mu+\ln 2) r^2\big(1+ O(\frac{1}{r})\big),$
			\item[ii)] $\mathbb E_{\mu}( |B_r|) = \frac{1}{3}(\mu+\ln 2) r^3\big(1+O(\frac{1}{r})\big)$.
		\end{itemize}
	\end{cor}
	\begin{proof}  It is clear that $ii)$ follows from $i)$ and \eqref{DB}. To establish i), we use that $\mathbb E_\rho(|D_r|) =1$ for all $r$, since $\rho$ is associated with a critical BGW process. By Proposition \ref{pushforward} ii) and the relation 
		$$
		D_r(T) = D_r(\chi(T)) + \sum_{(i,n):|i|<r} |D_{r-|i|+1}(T_{(i,n)})|\,, ~T \in \Omega_{\infty},
		$$
		with notation as above,  we get
		$$
		{\mathbb E}_\mu (|D_r|\,\mid \chi^{-1}(\cB^s_{\frac{1}{r}}(T^s))) = \sigma(B_r(T^s))+ |D_r(T^s)|\,, 
		$$
		where $\sigma$ is given by \eqref{sigmaform}. Integrating over $T^s$ then gives 
		$$
		\mathbb E_{\mu}(|D_r|) = \tilde{\mathbb E}_{\mu+\ln 2}(2 |B_r| -r)\,,
		$$
		and so relation i) follows from Corollary 4.9 ii) in \cite{meltem2021height}. 
	\end{proof}
	
	We conclude this section with stating the following almost sure result on the volume growth w.r.t. $\tau^{(\mu)}$, whose proof is basically identical to that of Theorem 4.11 in \cite{meltem2021height}, given Proposition \ref{pushforward}.
	
	\begin{prop}\label{asvolgrowth} For each $\mu>\mu_0$, there exist constants $C''_1, C''_2>0$ and for $\tau^{(\mu)}$-almost every $T\in\cT$ a number $r_0(T)\in\mathbb N$, such that 
		\bb\label{ashaus} \nonumber
		C''_1\cdot r^3\;\leq\; |B_r(T)|\;\leq C''_2\cdot r^3 \log r\,, \quad r\geq r_0(T)\,.
		\ee
		In particular, it holds that $d_h=3$ for $\tau^{(\mu)}$-a.e. tree $T$.
	\end{prop}

	\subsection*{Acknowledgements}
	
	The authors acknowledge support from Villum Fonden via the QMATH Centre of Excellence (Grant no.~10059).
	
	\section*{Appendix}
	
	\subsection*{Proof of Proposition \ref{prop:ZNln2}}	
	
	Defining 
	$$
	z=\sqrt{1-4g}
	$$
	and using \eqref{eq:recursionY} and \eqref{eq:sol1}, we get
	$$
	Y_m(g) = \prod_{l=1}^m X_l(g) = \frac{2^{-m}z}{\frac{1+z}{2}(1-z)^{-m} - \frac{1-z}{2}(1+z)^{-m}}\,.
	$$
	Hence, we have
	\bb\label{eq:Wa}
	W^{(\mu_0)}(g)  = \sum_{m=1}^\infty \frac{z}{D_m(z)}\,,
	\ee
	where 
	$$
	D_m(z) =  \frac{1+z}{2}(1-z)^{-m} - \frac{1-z}{2}(1+z)^{-m}\,.
	$$
	Note that $D_m$ is analytic in the unit disc. A straightforward power series expansion of the right-hand side gives 
	\bb\label{eq:powerDm}
	D_m(z) = (m+1) \cdot z \cdot \Big( 1+ \sum_{k \geq 1} b^{m} _{2k} (mz)^{2k} \Big)\,,
	\ee
	and that there exists a constant $K>0$ such that the coefficients $b^m _{2k}$ fulfill
	\bb\label{eq:b_2k1} \nonumber
	|b_{2k} ^m| \leq K\,,\quad k,m\geq 1 \,.
	\ee
	As observed in Section 2, $D_m$ is non-vanishing for $|z| < \tan\frac{\pi}{m+1}$ and by inverting the series \eqref{eq:powerDm} one obtains that for such values of $z$ it holds that 
	\bb \nonumber 
	\frac{z}{D_m(z)} = \frac{1}{m+1} \Big( 1+ \sum_{k \geq 1 } c_{2k}^m(mz)^{2k} \Big)\,,
	\ee
	where  the coefficients $c_{2k}^m$ fulfill
	\bb\label{eq:lem3b}
	|c_{2k} ^m| \leq (K+1)^k\, , \quad  k,m\geq 1 \,.
	\ee
	We refer to the proofs of Lemmas 4.1 and 4.2 in \cite{meltem2021height} for detailed arguments in an analogous situation.
	
	Consider first the contribution $S_1$ to the sum \eqref{eq:Wa} from $m\leq \frac{r}{|z|}$, for a given $r>0$ to be further specified later. With notation as above we have
	\bb \nonumber
	S_1 = \sum \limits_{m|z|\leq r}  \Big(\frac{1}{m+1} + \sum \limits_{k=1} ^\infty \frac{c^m _{2k}}{m+1}(mz)^{2k} \Big)\,.
	\ee
	Here, the first term inside round parenthesis is harmonic and yields the contribution 
	$$
	\ln\frac{r}{|z|}-1+\gamma + O\big(\frac{|z|}{r}\big)\,,
	$$
	while the bound \eqref{eq:lem3b} implies that the remaining contribution to $S_1$ bounded by  $O(|z| r) $, provided  $r$ fulfills
	\bb\label{eq:cond1}
	r<(K+1)^{-\frac 12}\,,
	\ee
	which we assume holds in the following. Thus, we obtain
	\bb \label{S11}
	S_1 =   \ln\frac{r}{|z|} -1+\gamma + O(r^2) +  O\big(\frac{|z|}{r}\big)\,.
	\ee
	The remaining contribution 
	\bb \nonumber
	S = \sum \limits_{m|z|> r}  \frac{z}{D_m(z)}\,.
	\ee
	to \eqref{eq:Wa} can be approximated by the integral
	$$
	I = \int \limits_{\frac{r}{|z|}} ^\infty \frac{z}{\sinh(sz)}ds\,,
	$$
	which can be rewritten as follows. First, note that  
	\bb \label{claimint}
	\int \limits_r ^\infty \frac{1}{\sinh t}  dt = \zeta - \ln r + O(r^2)
	\ee
	holds 	for small $r>0$, where $\zeta$ is a constant, since obviously the integral in \eqref{claimint} is convergent for $r>0$, and the small-$r$ behaviour  follows by applying the Laurent expansion of $\frac{1}{\sinh t}$ in the interval $[r,\frac{\pi}{2}]$.
	On the other hand, using Cauchy's theorem, we can rewrite the integral in \eqref{claimint} as a line integral along the circular arc $C_r$  of radius $r$ centered at $0$ connecting $r$ and $r\frac{z}{|z|}$ and along the half line $\ell_z: s \rightarrow sz$ with endpoint at $r\frac{z}{|z|}$ inside the wedge $V_a$, and get 
	\bb\label{cauchy2}
	\int \limits_{\frac{r}{|z|}}^\infty \frac{z}{\sinh(sz)}ds = \zeta - (\ln r +i\,\text{Arg}\,z) + O(r^2)\,,
	\ee
	where we have used that
	$$
	\int  \limits_{C_r} \frac{1}{\sinh w}dw  = i\,\text{Arg} z + O(r^2) \,.
	$$  
	Taking into account \eqref{cauchy2} in \eqref{S11}, we conclude that
	
	\begin{align}  \label{2piece2}
		\W^{(\mu_0)} (g) = \zeta + \gamma-1 - \ln z + (S - I)  +  O(r^2) +  O\big(\frac{|z|}{r}\big)\,.
	\end{align}
	
	Next, we proceed to estimate $|S-I|$ by splitting the summation and integration domains corresponding to $r < m|z|\leq \frac 1r$ and $m|z| > \frac 1r$ and similarly for $s$. Calling the corresponding sums and integrals $S_2, S_3$ and $I_2, I_3$, respectively, we first rewrite   
	$$
	S_2-I_2 = A + B \,,
	$$
	where 
	\begin{align*}
		A&:= \sum \limits_{r<m|z|\leq r^{-1}} \Big(\frac{z}{D_m(z)} - \frac{z}{\sinh(mz)} \Big) \\
		B&:= \sum \limits_{r<m|z|\leq r^{-1}} \frac{z}{\sinh(mz)} - \int \limits _{\frac{r}{|z|}} ^{\frac{1}{r|z|}} \frac{z ds}{\sinh(sz)} 
	\end{align*}
	In order to estimate $|A|$ and $|B|$ one may use that, for given $a>0$ and $0<\epsilon <1$, there exist constants $K_0, K_1,\delta_0>0$ and $m_0\in\mathbb N$ such that
	$$
	\sinh mx \leq |\sinh mz|  \leq K_0 \sinh mx\,, \qquad \frac 12\sinh(1-\epsilon)mx \leq |D_m(z)| \leq K_1\sinh(1+\epsilon)mx
	$$
	for all $z\in V_a$ such that $|z|<\delta_0$ and $m\geq m_0$. These inequalities are analogous to those given in Lemma 3.3 of \cite{meltem2021alpha} and their proofs are easily adapted. Following the proof of Theorem 3.4 of \cite{meltem2021alpha}, one obtains  the bounds 
	\bb \label{An}
	|A| \leq \text{cst} \cdot \frac{|z|}{r^3} \,, \quad |B| \leq \text{cst} \cdot \frac{|z|}{r^3} \,,
	\ee
	provided 
	\bb\label{eq:cond2}
	|z| < {\rm min}\{\delta_0,\frac{r}{m_0}\}\,,
	\ee		
	where $\delta_0$ here corresponds to, say, choosing $\epsilon=\frac 12$ above.
	
	For $S_3 - I_3$, on the other hand, one easily sees that there exist constants $c_1, c_2>0$ such that
	\bb \label{s3s3ni3}
	|S_3| \leq  \exp\big(-\frac{c_1}{r}\big)\,, \qquad  |I_3| \leq  \exp\big(-\frac{c_2}{r}\big)
	\ee
	for $0<r<1$.
	Collecting \eqref{2piece2}, \eqref{An} and \eqref{s3s3ni3} and defining 
	$$
	c_0 = \zeta + \gamma -1\,,
	$$
	we conclude that 
	\bb 
	\big| \W^{(\mu_0)} (g) - c_0 + \ln z\big|
	\,\leq\, \text{cst} \cdot \Big(r^2 + \frac{|z|}{r^{3}}\Big)\,, \label{snfinal}
	\ee
	provided $z$ and $r$  fulfill \eqref{eq:cond1} and \eqref{eq:cond2}.  
	Choosing $r=|z|^\beta$, where $0<\beta<1$, these conditions are satisfied for $|z|$ small enough, and setting $\beta=\frac 15$ 
	it follows from \eqref{snfinal} that \eqref{ZNln2} holds.
	
	This concludes the proof of Proposition \ref{prop:ZNln2}. \qed

	\bibliographystyle{abbrv}

\end{document}